\documentclass[11pt]{amsart}

\usepackage[T1]{fontenc}
\usepackage{lmodern}
\usepackage[margin=2.5cm]{geometry}
\usepackage{amssymb,mathtools}
\usepackage{enumitem}
\usepackage{tikz}
\usepackage{float}
\tikzset{
  gwire/.style={line width=.7pt,line cap=round,line join=round},
  gdash/.style={gwire,dash pattern=on 2.3pt off 1.6pt},
  gover/.style={preaction={draw=white,line width=4pt,solid,line cap=butt}},
  gbox/.style={draw,fill=white,line width=.7pt,inner sep=1.5pt,minimum height=.44cm,minimum width=.68cm,outer sep=0pt,font=\small},
  gdelta/.style={gbox,fill=black!10},
  glab/.style={font=\small,inner sep=1.3pt}
}
\usepackage[colorlinks=true,linkcolor=blue!55!black,citecolor=blue!55!black,urlcolor=blue!55!black]{hyperref}
\hypersetup{
  pdftitle={Localization of weakly integral objects in braided fusion categories},
  pdfauthor={Cesar Galindo}
}

\newtheorem{theorem}{Theorem}[section]
\newtheorem{proposition}[theorem]{Proposition}
\newtheorem{lemma}[theorem]{Lemma}
\newtheorem{corollary}[theorem]{Corollary}
\newtheorem{conjecture}[theorem]{Conjecture}
\theoremstyle{definition}
\newtheorem{definition}[theorem]{Definition}
\theoremstyle{remark}
\newtheorem{remark}[theorem]{Remark}
\numberwithin{equation}{section}
\newcommand{\C}{\mathbb{C}}
\newcommand{\Z}{\mathbb{Z}}
\newcommand{\cC}{\mathcal{C}}
\newcommand{\cD}{\mathcal{D}}
\newcommand{\cA}{\mathcal{A}}
\newcommand{\cE}{\mathcal{E}}
\newcommand{\cO}{\mathcal{O}}
\newcommand{\one}{\mathbf{1}}
\newcommand{\Vect}{\mathrm{Vec}}
\DeclareMathOperator{\Hom}{Hom}
\DeclareMathOperator{\End}{End}
\DeclareMathOperator{\FPdim}{FPdim}
\DeclareMathOperator{\Tr}{Tr}
\DeclareMathOperator{\tr}{tr}
\DeclareMathOperator{\RT}{RT}
\newcommand{\id}{\mathrm{id}}
\newcommand{\coev}{\mathrm{coev}}
\newcommand{\ad}{\mathrm{ad}}

\title[Localization of weakly integral objects]{Localization of weakly integral objects in braided fusion categories}
\author{C\'esar Galindo}
\address{Departamento de Matem\'aticas, Universidad de los Andes, Bogot\'a, Colombia}
\email{cn.galindo1116@uniandes.edu.co}
\date{November 2026}
\subjclass[2020]{18M20, 18M15, 20F36, 16T25, 57K16, 81P68}
\keywords{braided fusion category, braid group representation, Yang--Baxter operator, localization, Frobenius--Perron dimension, weakly integral, Turaev enhancement, topological quantum computation}

\AtBeginDocument{\lefthyphenmin=2 \righthyphenmin=3}

\begin{document}

\begin{abstract}
An object $X$ of a braided fusion category is called weakly integral if $\FPdim(X)^2$ is an integer. We prove that every nonzero weakly integral object admits a localization, that is, a Yang--Baxter operator $R$ on a finite-dimensional vector space $W$ such that, for every $n$, the algebra homomorphisms from the group algebra of the braid group $B_n$ to $\End(W^{\otimes n})$ and to $\End(X^{\otimes n})$, defined by $R$ and by the braiding, have the same kernel. The operator $R$ can be chosen unitary when the braided fusion category is unitary. This proves the `if' direction of a conjecture of Rowell and Wang, motivated by topological quantum computation, which states that a simple object admits a localization if and only if it is weakly integral. For a simple weakly integral object of a unitary braided fusion category, the operator $R$ that we construct has an enhancement in the sense of Turaev, and the resulting knot invariant is the Reshetikhin--Turaev invariant of knots colored by $X$.
\end{abstract}

\maketitle
\enlargethispage{3pt}

\section{Introduction}\label{sec:intro}

In topological quantum computation, information is stored in the state space of a collection of anyons and is processed by exchanging them \cite{Kit03,FKLW03,NSSFD08,RW18}. In the algebraic model, the types of anyons are the simple objects of a unitary braided fusion category $\cC$. The state space of $n$ anyons of type $X$ with total charge $Y$ is the space $\Hom_\cC(Y,X^{\otimes n})$, and the exchange of the anyons in positions $i$ and $i+1$ acts on it through the braiding $c_{X,X}$. These exchanges define a representation of the braid group $B_n$ on $\Hom_\cC(Y,X^{\otimes n})$, and its matrices are the gates of the computation.

The quantum circuit model is organized differently. Its state space is a tensor power $W^{\otimes n}$ of a fixed finite-dimensional Hilbert space $W$, and each elementary gate acts on a bounded number of factors. A braid group representation of this form, in which every generator acts by one and the same operator on two adjacent factors, is given by a Yang--Baxter operator. A \emph{Yang--Baxter operator} on a vector space $W$ is an invertible linear operator $R$ on $W\otimes W$ such that
\[
(R\otimes\id_W)(\id_W\otimes R)(R\otimes\id_W)=(\id_W\otimes R)(R\otimes\id_W)(\id_W\otimes R).
\]
It defines a representation $\rho^R_n$ of $B_n$ on $W^{\otimes n}$ by
\[
\rho^R_n(\sigma_i):=\id_W^{\otimes(i-1)}\otimes R\otimes\id_W^{\otimes(n-i-1)},\qquad 1\leq i\leq n-1 .
\]
We say that the representation $\rho^R_n$ is local, because each generator $\sigma_i$ acts only on two adjacent factors of $W^{\otimes n}$, and all the generators act by the same operator $R$. The representations that come from anyons are not local in any evident way, because the spaces $\Hom_\cC(Y,X^{\otimes n})$ are not tensor powers of a fixed vector space. Their dimensions $a_{n,Y}:=\dim\Hom_\cC(Y,X^{\otimes n})$ satisfy $\sum_Ya_{n,Y}\FPdim(Y)=d^n$, where the sum is over the simple objects $Y$ and $d:=\FPdim(X)$ is the Frobenius--Perron dimension of $X$, that is, the quantum dimension of the anyon. The number $d$ is rarely an integer. For the Fibonacci anyon, for example, $d$ is the golden ratio and the numbers $a_{n,Y}$ are Fibonacci numbers.

Rowell and Wang asked when the braid group representations of an anyon $X$ can nevertheless be realized by a Yang--Baxter operator \cite{RW12}. The representation $\rho^R_n$ on $W^{\otimes n}$ and the representation of $B_n$ on the spaces $\Hom_\cC(Y,X^{\otimes n})$ cannot be isomorphic in general, because their dimensions are different. The requirement is weaker. For every $n$, the representation $\rho^R_n$ on $W^{\otimes n}$ must contain the same irreducible representations of $B_n$ as the direct sum of the spaces $\Hom_\cC(Y,X^{\otimes n})$, where $Y$ runs over the simple objects, possibly with different multiplicities, and the same operator $R$ must be used for every $n$. In that case every irreducible representation of $B_n$ that occurs in a space $\Hom_\cC(Y,X^{\otimes n})$ is realized inside $W^{\otimes n}$, where an elementary exchange is implemented by applying the same gate $R$ to two adjacent tensor factors \cite[Definition~5.2]{RW18}. The simulation of topological quantum field theories by quantum circuits in general is treated in \cite{FKW02}.

The precise notion is the following. Throughout the paper, fusion categories are over $\C$. Let $X$ be an object of a braided fusion category $\cC$, and let
\[
\rho^X_n\colon\C B_n\longrightarrow\End_\cC(X^{\otimes n}),\qquad \sigma_i\longmapsto\id_X^{\otimes(i-1)}\otimes c_{X,X}\otimes\id_X^{\otimes(n-i-1)}
\]
be the algebra homomorphism defined by the braiding. A \emph{localization} of $X$ is a pair $(W,R)$, where $W$ is a nonzero finite-dimensional vector space and $R$ is a Yang--Baxter operator on $W$, such that for every $n\geq2$ there is an injective algebra homomorphism $\phi_n$ from the image of $\rho^X_n$ to $\End(W^{\otimes n})$ with $\phi_n\circ\rho^X_n=\rho^R_n$ \cite[Definition~2.3]{RW12}, \cite[Definition~1.2]{GHR13}. Equivalently, $\rho^X_n$ and $\rho^R_n$ have the same kernel in the group algebra $\C B_n$ for every $n\geq2$. In the unitary case, where these representations are completely reducible, this means that $W^{\otimes n}$ and the direct sum of the spaces $\Hom_\cC(Y,X^{\otimes n})$ contain the same irreducible representations of $B_n$. A localization is \emph{unitary} when $W$ is a Hilbert space and $R$ is unitary. An object that admits a localization is called \emph{localizable}.

\begin{conjecture}[Rowell--Wang]\label{conj:RW-intro}
Let $X$ be a simple object of a braided fusion category $\cC$.
\begin{enumerate}[label=\textup{(\alph*)}]
\item The object $X$ admits a localization if and only if $\FPdim(X)^2$ is an integer \textup{\cite[Conjecture~4.1]{RW12}}.
\item If $\cC$ is unitary, then $X$ admits a unitary localization if and only if $\FPdim(X)^2$ is an integer \textup{\cite[Conjecture~6.3]{RW18}}.
\end{enumerate}
\end{conjecture}

We say that an object $X$ of a fusion category is \emph{weakly integral} if $\FPdim(X)^2$ is an integer. A simple object of a braided fusion category is weakly integral precisely when the fusion subcategory that it generates is weakly integral, that is, its Frobenius--Perron dimension is an integer (Lemma~\ref{lem:integrality}). In these terms, the conjecture says that the localizable anyons are the weakly integral ones.

Part~(a) of Conjecture~\ref{conj:RW-intro} is stated in \cite{RW12} for braided fusion categories that need not be unitary, and we use its algebraic formulation from \cite[Conjecture~1.4]{GHR13}. For part~(b) of Conjecture~\ref{conj:RW-intro} see also \cite[Conjecture~1.2]{GR26}. Two examples illustrate the statement. The Ising anyon $\sigma$, with $d=\sqrt2$, admits a unitary localization of dimension $2$ \cite{FRW06,RW12}. The Fibonacci anyon $\tau$, for which $d^2$ is not an integer, admits no localization. Indeed, its braid group images generate the algebras $\End(\tau^{\otimes n})$ \cite[Section~6.2]{DRW16}, and the claim follows from \cite[Theorem~4.5]{RW12}.

The conjecture is part of a description of the anyons whose braiding gates alone are not universal for quantum computation. For the Jones representations, that is, the braid group representations of the generating simple object of the braided fusion categories $SU(2)_k$, Rowell and Wang showed that a unitary localization exists only at the levels $k=1$, $2$ and $4$, which are also the only levels at which braiding alone is not universal \cite{RW12}. In general, the Property~F conjecture predicts that the braid group images $\rho^X_n(B_n)$ of a simple object $X$ are finite for every $n\geq2$ if and only if $X$ is weakly integral \cite{NR11}, \cite[Conjecture~6.2]{RW18}. Rowell and Wang also conjectured that a unitary Yang--Baxter operator of finite order has finite braid group images, modulo scalars \cite[Conjecture~3.1]{RW12}, and they proposed Conjecture~\ref{conj:RW-intro} as the link between these two statements. For the Property~F conjecture, finiteness is known when the fusion subcategory generated by $X$ is weakly group-theoretical \cite[Theorem~7.2]{GN21}. In the other direction, a simple object with finite braid group images is weakly integral when the braided fusion category is pseudo-unitary \cite[Corollary~1.2]{LR26}, and when the fusion subcategory generated by $X$ is nondegenerate and pivotal or has a pseudo-unitary Galois conjugate \cite[Corollary~5.6]{ABCPT26}. These results concern finite images. It remains open whether every localizable simple object is weakly integral.

Two general results were known. Rowell and Wang proved the `only if' direction when the image of the braid group generates the algebra $\End_\cC(X^{\otimes n})$ for every $n\geq2$ \cite[Theorem~4.5]{RW12}; see also \cite[Theorem~6.3]{GHR13}. In the other direction, a finite-dimensional representation of a semisimple quasi-triangular Hopf algebra is localized by its braiding, regarded as an operator on the tensor square of the underlying vector space \cite[Proposition~4.17]{GHR13}, and in this case the Frobenius--Perron dimension of the object is an integer. For objects of nonintegral dimension, localizations were known only in specific families \cite{FRW06,RW12,GR26}. In general only weaker notions had been realized, such as the weak localizations of \cite[Definition~4.20 and Proposition~4.21]{GHR13}, in which the tensor powers of a vector space are replaced by tensor powers of a bimodule over a semisimple algebra.

\subsection*{The localization theorem}
We prove the `if' direction of both parts of Conjecture~\ref{conj:RW-intro}. The sufficiency statement holds for arbitrary nonzero objects. The converse in Conjecture~\ref{conj:RW-intro} remains open for simple objects, and we do not address~it. For objects that are not simple, the condition is not necessary. Indeed, if $(W,R)$ is a localization of an object $Y$, then $\one\oplus Y$ is localized by $\C\oplus W$ with the operator that acts by $R$ on $W\otimes W$ and by the flip on the other three summands of $(\C\oplus W)^{\otimes2}$, which is a Yang--Baxter operator \cite[Proposition~3.8]{GR26}. The reason is that both tensor powers decompose into summands indexed by the sets $S$ of positions occupied by $Y$ and by $W$, and in both cases a braid $b$ sends the summand of $S$ to the summand of the image of $S$ under the permutation of $b$, through the braid obtained from $b$ by deleting the other strands. For the Ising anyon this gives a unitary localization of $\one\oplus\sigma$ of dimension $3$, and $\FPdim(\one\oplus\sigma)^2=3+2\sqrt2$ is not an integer.

\begin{theorem}\label{thm:rw}
Let $X$ be a nonzero weakly integral object of a braided fusion category $\cC$. Then $X$ admits a localization. If $\cC$ is unitary, then $X$ admits a unitary localization.
\end{theorem}

\begin{corollary}\label{cor:rw-intro}
Let $\cC$ be a braided fusion category.
\begin{enumerate}[label=\textup{(\alph*)}]
\item If $\cC$ is unitary and weakly integral, then every simple object of $\cC$ admits a unitary localization.
\item If $X$ is a simple object of $\cC$ and the fusion subcategory generated by $X$ is weakly group-theoretical \cite[Definition~9.8.1]{EGNO15}, then $X$ admits a unitary localization. In particular, if $\cC$ is weakly group-theoretical, then every simple object of $\cC$ admits a unitary localization.
\end{enumerate}
\end{corollary}

In part~(b) of Corollary~\ref{cor:rw-intro}, the braided fusion category $\cC$ is not assumed to be unitary. More generally, a nonzero weakly integral object $X$ admits a unitary localization whenever the fusion subcategory generated by $X\otimes X^*$ is tensor equivalent to a unitary fusion category (Corollary~\ref{cor:unitarizable}).

The operator $R$ that we construct has finite order (Remark~\ref{rem:finite-order}). Hence, for unitary $\cC$, the conjecture of Rowell and Wang on unitary Yang--Baxter operators of finite order would now imply that the groups $\rho^X_n(B_n)$ are finite for every weakly integral object $X$, as the Property~F conjecture predicts. Together with \cite[Corollary~1.2]{LR26}, it would give the Property~F conjecture for the simple objects of unitary braided fusion categories.

\subsection*{The main construction}
The image of $\rho^X_n$ is in general a proper subalgebra of the endomorphism algebra $\End_\cC(X^{\otimes n})$. We prove Theorem~\ref{thm:rw} by representing the whole endomorphism algebra on $W^{\otimes n}$, for every $n$, in such a way that the tensor product of two endomorphisms is sent to the tensor product of the corresponding operators. Theorem~\ref{thm:main} is the precise statement.

\begin{theorem}\label{thm:main}
Let $X$ be a nonzero object of a braided fusion category $\cC$. Put $d:=\FPdim(X)$ and $E_n:=\End_\cC(X^{\otimes n})$, and let $\cA$ be the fusion subcategory of $\cC$ generated by $X\otimes X^*$. Suppose that $d^2$ is an integer. Then $m:=d^2\FPdim(\cA)$ is an integer, and there are a vector space $W$ of dimension $m$ and injective unital algebra homomorphisms $\Phi_n\colon E_n\to\End(W^{\otimes n})$, $n\geq1$, with the following properties.
\begin{enumerate}[label=\textup{(\alph*)}]
\item For all $n,k\geq1$, $f\in E_n$ and $g\in E_k$,
\begin{equation}\label{eq:tensor-compatible}
\Phi_{n+k}(f\otimes g)=\Phi_n(f)\otimes\Phi_k(g).
\end{equation}
\item Regard $W^{\otimes n}$ as an $E_n$-module via $\Phi_n$. For each simple constituent $Y$ of $X^{\otimes n}$, the simple module $\Hom_\cC(Y,X^{\otimes n})$ occurs in $W^{\otimes n}$ with multiplicity $(m/d)^n\FPdim(Y)$. Consequently,
\begin{equation}\label{eq:trace-identity}
m^{-n}\tr\bigl(\Phi_n(f)\bigr)=d^{-n}\Tr(f),\qquad f\in E_n,
\end{equation}
where $\tr$ is the trace of linear maps and $\Tr$ is the trace of the canonical spherical structure of the fusion subcategory generated by $X$ \textup{(}Subsection~\textup{\ref{subsec:spherical})}.
\item If $\cC$ is unitary, then $W$ has an inner product and the homomorphisms $\Phi_n$ can be chosen in such a way that $\Phi_n(f^*)=\Phi_n(f)^*$ for all $n\geq1$ and $f\in E_n$.
\end{enumerate}
\end{theorem}

Theorem~\ref{thm:rw} follows by taking $R:=\Phi_2(c_{X,X})$. By \eqref{eq:tensor-compatible}, we have
\[
\Phi_n\bigl(\rho^X_n(\sigma_i)\bigr)=\id_W^{\otimes(i-1)}\otimes R\otimes\id_W^{\otimes(n-i-1)}=\rho^R_n(\sigma_i),\qquad 1\leq i\leq n-1.
\] The case $n=3$ gives the Yang--Baxter equation, and the restriction of the injective homomorphism $\Phi_n$ to the image of $\rho^X_n$ is the homomorphism required in the definition of a localization (Lemma~\ref{lem:localization}). Part~(b) of Theorem~\ref{thm:main} describes $W^{\otimes n}$ as a module over $E_n$. It is the direct sum of the spaces $\Hom_\cC(Y,X^{\otimes n})$, where $Y$ runs over the simple constituents of $X^{\otimes n}$ and the space $\Hom_\cC(Y,X^{\otimes n})$ is repeated $(m/d)^n\FPdim(Y)$ times. Restricting this decomposition along $\rho^X_n$ gives a decomposition of $W^{\otimes n}$ as a representation of $B_n$. For two different simple objects $Y$ and $Y'$, the spaces $\Hom_\cC(Y,X^{\otimes n})$ and $\Hom_\cC(Y',X^{\otimes n})$ can have common irreducible constituents as representations of $B_n$, and then their multiplicities add.

If the fusion subcategory generated by $X$ had a fiber functor, that is, a tensor functor to vector spaces, then this functor would give homomorphisms $\Phi_n$ with $\dim W=d$, as in \cite[Proposition~4.17]{GHR13}. A fiber functor exists only if $d$ is an integer. The homomorphisms of Theorem~\ref{thm:main} are defined only on the endomorphisms of the tensor powers of $X$, and they need not come from a tensor functor.

The construction has two ingredients. The first is to pair each copy of $X$ with a copy of its dual. The object $X\otimes X^*$ has dimension $d^2$, and when $d^2$ is an integer the fusion subcategory $\cA$ generated by $X\otimes X^*$ is integral, that is, the Frobenius--Perron dimension of each of its objects is an integer. The braiding allows us to bring the copies of $X$ in $(X\otimes X^*)^{\otimes n}$ together, apply $f\in E_n$, and return them to their places. This gives homomorphisms $E_n\to\End_\cA\bigl((X\otimes X^*)^{\otimes n}\bigr)$ that are compatible with tensor products (Section~\ref{sec:sufficiency}).

The second ingredient is the regular object of the integral fusion category $\cA$,
\[
R_\cA:=\bigoplus_{Y}\FPdim(Y)\,Y,
\]
where $Y$ runs over the simple objects of $\cA$ and $\FPdim(Y)\,Y$ denotes the direct sum of $\FPdim(Y)$ copies of $Y$. It satisfies
\[
V\otimes R_\cA\cong\FPdim(V)\,R_\cA
\]
for every object $V$ of $\cA$. The copies of $R_\cA$ in $V\otimes R_\cA$ are parametrized by the space $\Hom_\cA(R_\cA,V)$, whose dimension is $\FPdim(V)$, and a single morphism $\delta\colon R_\cA\to R_\cA\otimes R_\cA$ gives this parametrization for every $V$ at once (Lemma~\ref{lem:beta}). This choice produces homomorphisms
\[
\End_\cA\bigl((X\otimes X^*)^{\otimes n}\bigr)\longrightarrow\End(W^{\otimes n}),\qquad W:=\Hom_\cA(R_\cA,X\otimes X^*\otimes R_\cA),
\]
that are again compatible with tensor products (Section~\ref{sec:regular}). The homomorphisms $\Phi_n$ are the composites. For part~(b) of Theorem~\ref{thm:main}, the construction identifies $W^{\otimes n}$ with a space of morphisms into $X^{\otimes n}\otimes(X^*)^{\otimes n}\otimes R_\cA^{\otimes n}$, on which $E_n$ acts through the first factor. Decomposing $X^{\otimes n}$ into simple objects then exhibits the space that carries each multiplicity, and its dimension is a product of Frobenius--Perron dimensions.

\subsection*{Knot invariants}
One of the questions that led Rowell and Wang to localizations is whether a unitary Yang--Baxter operator can produce, through a weighted trace, the values of the Jones polynomial at roots of unity, since such an operator gives a direct way to approximate these values in the quantum circuit model \cite{RW12}. Turaev showed that a Yang--Baxter operator $R$ gives a link invariant when it has an enhancement, that is, an operator $\mu$ on $W$ and two scalars for which the traces $\tr(\rho^R_n(b)\mu^{\otimes n})$, suitably normalized, are invariant under the Markov moves \cite{Tur88}. For a simple weakly integral object $X$ of a unitary braided fusion category, the operator $R=\Phi_2(c_{X,X})$ obtained from Theorem~\ref{thm:main} has an enhancement with a scalar operator $\mu$, and the resulting invariant is the Reshetikhin--Turaev invariant of knots colored by $X$.

For a Hilbert space $W$ of dimension $m$ we denote by $\tau_n:=m^{-n}\tr$ the normalized trace on $\End(W^{\otimes n})$, and by $\id\otimes\tau_1$ and $\tau_1\otimes\id$ the linear maps $\End(W\otimes W)\to\End(W)$ that send $x\otimes y$ to $\tau_1(y)\,x$ and to $\tau_1(x)\,y$, respectively. In these expressions $\id$ is the identity map of $\End(W)$.

\begin{theorem}\label{thm:enhanced-intro}
Let $X$ be a simple object of a unitary braided fusion category $\cC$, put $d:=\FPdim(X)$, and suppose that $d^2$ is an integer. Let $\theta_X$ be the twist of $X$ for the canonical ribbon structure of the fusion subcategory generated by $X$ \textup{(}Subsection~\textup{\ref{subsec:ribbon})}. Let $W$ be the Hilbert space of dimension $m$ and $\Phi_n$ the homomorphisms of Theorem~\textup{\ref{thm:main}(c)}, and put $R:=\Phi_2(c_{X,X})$.
\begin{enumerate}[label=\textup{(\alph*)}]
\item The normalized partial traces of $R^{\pm1}$ are the scalar operators
\[
(\id\otimes\tau_1)(R^{\pm1})=(\tau_1\otimes\id)(R^{\pm1})=\theta_X^{\pm1}d^{-1}\id_W .
\]
In particular, $\tau_2(R)=\theta_X/d$, and the operator $R$ determines $d=|\tau_2(R)|^{-1}$ and $\theta_X=d\,\tau_2(R)$.
\item The quadruple $(R,(d/m)\id_W,\theta_X,1)$ is an enhanced Yang--Baxter operator in the sense of Turaev \textup{\cite{Tur88}}. If $b\in B_n$ is a braid whose closure $\widehat b$ is a knot, then
\begin{equation}\label{eq:local-rt}
\RT_X(\widehat b)=\theta_X^{-w(b)}\Bigl(\frac{d}{m}\Bigr)^n\tr\bigl(\rho^R_n(b)\bigr),
\end{equation}
where $w(b)$ is the exponent sum of $b$ and $\RT_X(\widehat b)$ is the Reshetikhin--Turaev invariant of the zero-framed knot $\widehat b$ colored by $X$, with value $d$ on the unknot.
\end{enumerate}
\end{theorem}

The proof uses only the properties \eqref{eq:tensor-compatible}, \eqref{eq:trace-identity} and $\Phi_n(f^*)=\Phi_n(f)^*$ of the homomorphisms $\Phi_n$. The knot formula follows from the trace identity \eqref{eq:trace-identity}. For part~(a) of Theorem~\ref{thm:enhanced-intro} we use that adding a kink to the closure of a braid multiplies the normalized trace by $\theta_X^{\pm1}/d$ (Figure~\ref{fig:stab}). This is the Markov property of the normalized traces of the representations $\rho^R_n$, and by a criterion of Lechner it implies that the partial trace of $R$ over the first factor is the scalar operator $\theta_Xd^{-1}\id_W$ \cite[Proposition~2.3]{Lec26}. The two partial traces of a unitary Yang--Baxter operator coincide, by a theorem of Conti and Lechner \cite[Theorem~5.10(a)]{CL21}.

The paper is organized as follows. Section~\ref{sec:preliminaries} recalls the notions that we use. Section~\ref{sec:regular} contains the construction with the regular object of an integral fusion category, and Section~\ref{sec:sufficiency} pairs each copy of $X$ with a copy of $X^*$ and proves parts (a) and (b) of Theorem~\ref{thm:main}. Section~\ref{sec:unitary} treats the unitary case and completes the proofs of Theorems~\ref{thm:rw} and~\ref{thm:main} and of Corollary~\ref{cor:rw-intro}. Section~\ref{sec:turaev} proves Theorem~\ref{thm:enhanced-intro}. For the articles that are available on the arXiv, the numbers of the cited statements are those of the arXiv versions indicated in the bibliography.

%======================================================================
\section{Preliminaries}\label{sec:preliminaries}

This section fixes the notation and recalls the facts about fusion categories that are used in the paper. The statements used repeatedly are the isomorphisms \eqref{eq:regular-absorbs} for the regular object, the description \eqref{eq:wedderburn} of the endomorphism algebra of an object and the trace formula (Lemma~\ref{lem:trace-formula}). Subsection~\ref{subsec:localization} contains the definition of a localization and Lemma~\ref{lem:localization}, which reduces Theorem~\ref{thm:rw} to the construction of the homomorphisms $\Phi_n$ of Theorem~\ref{thm:main}.

We work over $\C$ and write $\Vect$ for the category of finite-dimensional complex vector spaces. We write tensor products without parentheses, using the coherence isomorphisms to identify different bracketings \cite[Theorems~2.8.5 and~2.9.2]{EGNO15}. In a unitary fusion category these identifications are unitary. In particular, conjugating an endomorphism by a change of bracketing preserves adjoints. For every object $V$ we put $V^{\otimes0}:=\one$.

\subsection{Fusion categories}\label{subsec:fusion}

Let $\cC$ be a fusion category, that is, a finite semisimple rigid $\C$-linear tensor category with simple unit. We write $\cO(\cC)$ for a set of representatives of the isomorphism classes of its simple objects. For $Y\in\cO(\cC)$, the multiplicity of $Y$ in an object $V$ is $[V:Y]:=\dim\Hom_\cC(Y,V)$. We call $Y$ a \emph{constituent} of $V$ when this multiplicity is positive. We write $V^*$ for the left dual of $V$, with coevaluation $\coev_V\colon\one\to V\otimes V^*$ and evaluation $V^*\otimes V\to\one$, and ${}^*V$ for a right dual of $V$, with evaluation $V\otimes{}^*V\to\one$ \cite[Definitions~2.10.1 and~2.10.2]{EGNO15}.

A fusion subcategory is a full subcategory that is closed under isomorphisms, direct summands, finite direct sums, tensor products and duals, and contains the unit. The fusion subcategory generated by a collection of objects is the smallest fusion subcategory that contains them. For a single object $V$ we denote it by $\langle V\rangle$.

\begin{lemma}\label{lem:generated}
Let $V$ be a nonzero object of a fusion category. A simple object belongs to $\langle V\rangle$ if and only if it is a direct summand of $V^{\otimes n}$ for some $n\geq1$.
\end{lemma}

\begin{proof}
Let $\mathcal K$ be the full subcategory consisting of the objects all of whose simple constituents are direct summands of the powers $V^{\otimes n}$, $n\geq1$. It contains $V$, and it is closed under isomorphisms, direct sums and direct summands. It is also closed under tensor products. Indeed, a simple constituent of $A\otimes B$ is a constituent of $Y\otimes Y'$ for some simple constituents $Y$ of $A$ and $Y'$ of $B$, and if $Y$ and $Y'$ are direct summands of $V^{\otimes a}$ and $V^{\otimes b}$, then $Y\otimes Y'$ is a direct summand of $V^{\otimes(a+b)}$. By \cite[Corollary~4.11.4]{EGNO15}, $\mathcal K$ contains $\one$ and is closed under duals. Hence $\mathcal K$ is a fusion subcategory that contains $V$, and therefore $\langle V\rangle\subseteq\mathcal K$. Conversely, a fusion subcategory that contains $V$ contains the powers $V^{\otimes n}$ and their direct summands, and therefore it contains $\mathcal K$. Thus $\mathcal K=\langle V\rangle$.
\end{proof}

We write $d_V:=\FPdim(V)$ for the Frobenius--Perron dimension. It is additive on direct sums and multiplicative on tensor products, and $d_{V^*}=d_V$. For a simple object $Y$, the number $d_Y$ is a positive algebraic integer. A fusion category is \emph{integral} if $d_Y$ is an integer for every simple object $Y$, and \emph{weakly integral} if $\FPdim(\cC):=\sum_{Y\in\cO(\cC)}d_Y^2$ is an integer. Tensor functors between fusion categories preserve Frobenius--Perron dimensions, and hence $d_V$ can be computed in any fusion subcategory containing $V$ \cite[Sections~3.3 and~4.5 and Definition~9.6.1]{EGNO15}.

When $\cC$ is integral, the \emph{regular object} is $R_\cC:=\bigoplus_{Y\in\cO(\cC)} d_Y\,Y$, the direct sum of $d_Y$ copies of each simple object $Y$. Its dimension is $\FPdim(\cC)$ \cite[Definitions~6.1.6 and~6.1.7]{EGNO15}. In the Grothendieck ring of $\cC$ the class of the regular object satisfies $[V]\,[R_\cC]=d_V\,[R_\cC]$ for every object $V$ \cite[Proposition~6.1.11]{EGNO15}. In a fusion category two objects with the same class in the Grothendieck ring are isomorphic. Hence there are isomorphisms
\begin{equation}\label{eq:regular-absorbs}
V\otimes R_\cC\cong d_V\,R_\cC\qquad(V\in\cC),
\end{equation}
where $d_V\,R_\cC$ denotes the direct sum of $d_V$ copies of $R_\cC$; see also \cite[Proposition~3.15]{MT08}.

If $V$ is a nonzero object of $\cC$, then $\one$ is a direct summand of $V\otimes V^*$, because $\coev_V\colon\one\to V\otimes V^*$ is nonzero by rigidity and $\cC$ is semisimple. Choose $p\colon V\otimes V^*\to\one$ with $p\,\coev_V=\id_\one$. Then every morphism $f\colon A\to B$ is recovered from $f\otimes\id_V$ by the formula
\begin{equation}\label{eq:faithful-tensor}
f=(\id_B\otimes p)(f\otimes\id_V\otimes\id_{V^*})(\id_A\otimes\coev_V),
\end{equation}
and hence the functor $-\otimes V$ is faithful.

Finally, we recall the structure of the endomorphism algebra of an object $V$. For $f\in\End_\cC(V)$ and a simple object $Y$, let $f_Y$ be the linear map $v\mapsto fv$ on $\Hom_\cC(Y,V)$. Since $\cC$ is semisimple, the map
\begin{equation}\label{eq:wedderburn}
\End_\cC(V)\longrightarrow\bigoplus_Y\End\bigl(\Hom_\cC(Y,V)\bigr),\qquad f\longmapsto(f_Y)_Y,
\end{equation}
where $Y$ runs over the simple constituents of $V$, is an algebra isomorphism (cf.\ \cite[Lemma~II.4.2.2]{Tur10}). In particular, $\End_\cC(V)$ is semisimple, and its simple modules are the spaces $\Hom_\cC(Y,V)$.

For a finite-dimensional vector space $L$ and an object $Y$, we write $L\odot Y$ for the direct sum of copies of $Y$ with multiplicity space $L$. Thus every $\ell\in L$ gives a morphism $\iota_\ell\colon Y\to L\odot Y$ that depends linearly on $\ell$, and a morphism from $L\odot Y$ to an object $Z$ is the same as a linear map from $L$ to $\Hom_\cC(Y,Z)$. For every object $V$, the morphism
\begin{equation}\label{eq:canonical-decomposition}
\bigoplus_{Y\in\cO(\cC)}\Hom_\cC(Y,V)\odot Y\longrightarrow V
\end{equation}
that corresponds to the identity maps of the spaces $\Hom_\cC(Y,V)$ is an isomorphism, natural in $V$.

\subsection{Traces and pseudo-unitary fusion categories}\label{subsec:spherical}

A pivotal structure is a monoidal natural isomorphism $\mathfrak a\colon\id_\cC\to(-)^{**}$. We write $\Tr(f)$ for the resulting left trace of an endomorphism $f$, and $\dim(V):=\Tr(\id_V)$. The structure is \emph{spherical} if the left and right traces agree. A fusion category $\cC$ is \emph{pseudo-unitary} if its categorical dimension $\dim(\cC)$ \cite[Definition~7.21.3]{EGNO15} equals $\FPdim(\cC)$ \cite[Definition~9.4.4]{EGNO15}. It then has a unique spherical structure with $\dim(Y)=d_Y$ for every simple object $Y$, called the \emph{canonical spherical structure}. Every weakly integral fusion category is pseudo-unitary \cite[Sections~4.7 and~9.4 and Propositions~9.5.1 and~9.6.5]{EGNO15}.

We distinguish the categorical trace $\Tr$ from the trace $\tr$ of linear maps.

\begin{lemma}\label{lem:trace-formula}
In a spherical fusion category $\cC$, for every object $V$ and $f\in\End_\cC(V)$,
\[
\Tr(f)=\sum_{Y\in\cO(\cC)}\dim(Y)\tr(f_Y).
\]
\end{lemma}

\begin{proof}
The trace is a linear form on $\End_\cC(V)$ with $\Tr(gh)=\Tr(hg)$. By \eqref{eq:wedderburn}, the algebra $\End_\cC(V)$ is the product of the matrix algebras $\End(\Hom_\cC(Y,V))$, and a linear form with this property on a matrix algebra is a scalar multiple of the trace $\tr$. Hence there are scalars $c_Y$ such that $\Tr(f)=\sum_Yc_Y\tr(f_Y)$ for every $f$. To compute $c_Y$, let $u\colon Y\to V$ be a nonzero morphism, and let $p\colon V\to Y$ be a morphism with $pu=\id_Y$, which exists because $\cC$ is semisimple. For $e:=up$, the map $e_Y$ is an idempotent of rank one and $e_{Y'}=0$ for $Y'\neq Y$. Therefore $c_Y=\Tr(up)=\Tr(pu)=\dim(Y)$.
\end{proof}

\subsection{Braided fusion categories and braid groups}\label{subsec:braided}

We write $c_{A,B}\colon A\otimes B\to B\otimes A$ for the braiding. It is natural in both variables, and it satisfies the hexagon identities
\begin{equation}\label{eq:hexagons}
\begin{aligned}
c_{A,B\otimes C}&=(\id_B\otimes c_{A,C})(c_{A,B}\otimes\id_C),\\
c_{A\otimes B,C}&=(c_{A,C}\otimes\id_B)(\id_A\otimes c_{B,C}).
\end{aligned}
\end{equation}
The braid group $B_n$ is generated by $\sigma_1,\dots,\sigma_{n-1}$ with the relations $\sigma_i\sigma_{i+1}\sigma_i=\sigma_{i+1}\sigma_i\sigma_{i+1}$ and $\sigma_i\sigma_j=\sigma_j\sigma_i$ for $|i-j|\geq2$. For an object $X$ of a braided fusion category $\cC$ we write $E_n:=\End_\cC(X^{\otimes n})$. The braiding defines an algebra homomorphism
\[
\rho^X_n\colon\C B_n\to E_n,\qquad
\rho^X_n(\sigma_i)=\id_X^{\otimes(i-1)}\otimes c_{X,X}\otimes\id_X^{\otimes(n-i-1)} .
\]

We shall use the following fact about monoidal functors twice. Let $G$ be a monoidal functor between monoidal categories, with structure isomorphisms $\gamma_{A,B}\colon G(A)\otimes G(B)\to G(A\otimes B)$ \cite[Definition~2.4.1]{EGNO15}. We do not write the associativity constraints, and then the associativity axiom of $\gamma$ is the identity $\gamma_{A\otimes B,C}\,(\gamma_{A,B}\otimes\id_{G(C)})=\gamma_{A,B\otimes C}\,(\id_{G(A)}\otimes\gamma_{B,C})$. For an object $V$ put $\gamma_{n,k}:=\gamma_{V^{\otimes n},V^{\otimes k}}$, and define the iterated structure isomorphisms $\gamma_n\colon G(V)^{\otimes n}\to G(V^{\otimes n})$ by $\gamma_1:=\id_{G(V)}$ and $\gamma_{n+1}:=\gamma_{n,1}\,(\gamma_n\otimes\id_{G(V)})$. Then, for all $n,k\geq1$,
\begin{equation}\label{eq:gamma-iterated}
\gamma_{n+k}=\gamma_{n,k}\,\bigl(\gamma_n\otimes\gamma_k\bigr).
\end{equation}
Indeed, for $k=1$ this is the definition of $\gamma_{n+1}$, and if \eqref{eq:gamma-iterated} holds for $k$, then the associativity axiom gives
\[
\gamma_{n+k+1}=\gamma_{n+k,1}\,\bigl(\gamma_{n,k}\otimes\id\bigr)\bigl(\gamma_n\otimes\gamma_k\otimes\id\bigr)
=\gamma_{n,k+1}\,\bigl(\id\otimes\gamma_{k,1}\bigr)\bigl(\gamma_n\otimes\gamma_k\otimes\id\bigr)
=\gamma_{n,k+1}\,\bigl(\gamma_n\otimes\gamma_{k+1}\bigr).
\]
\subsection{Twists and Reshetikhin--Turaev invariants}\label{subsec:ribbon}

A spherical structure $\mathfrak a$ on a braided fusion category determines the ribbon twist $\theta$ by $\mathfrak a_V=u_V\theta_V$, where $u$ is the Drinfeld isomorphism \cite[Sections~8.9 and~8.10]{EGNO15}. In a pseudo-unitary braided fusion category we take the twist associated with its canonical spherical structure. For a simple object $X$ the automorphism $\theta_X$ is a scalar multiple of $\id_X$, and we also denote this scalar by $\theta_X$. It is a root of unity \cite[Corollary~8.18.2(i)]{EGNO15}.

Let $X$ be a simple object of a pseudo-unitary braided fusion category, let $b\in B_n$ be a braid and let $\widehat b$ be its closure, with the orientation given by the braid. The trace $\Tr(\rho^X_n(b))$ is the Reshetikhin--Turaev invariant of the link $\widehat b$, with its blackboard framing and with all components colored by $X$ \cite{RT90}, \cite[Theorem~I.2.5 and Corollary~I.2.7.2]{Tur10}. In particular, it depends only on the framed link. A positive twist of the framing of a component multiplies the invariant by $\theta_X$, and a negative twist multiplies it by $\theta_X^{-1}$ \cite[Theorem~I.2.5]{Tur10}. When $\widehat b$ is a knot, its blackboard framing is the exponent sum $w(b)$ of the braid. The invariant of the zero-framed knot is therefore
\begin{equation}\label{eq:rt}
\RT_X(\widehat b)=\theta_X^{-w(b)}\Tr(\rho^X_n(b)).
\end{equation}
The invariant of the zero-framed unknot is $\dim(X)$, since the unknot is the closure of the trivial braid on one strand. We shall use the following consequence. The closures of $\sigma_1$ and $\sigma_1^{-1}$ in $B_2$ are unknots with framings $1$ and $-1$, and hence \eqref{eq:rt} gives
\begin{equation}\label{eq:twist-trace}
\Tr(c_{X,X}^{\pm1})=\theta_X^{\pm1}\dim(X).
\end{equation}
For the negative sign this is \cite[Exercise~8.10.15]{EGNO15}.

\subsection{Unitary fusion categories}\label{subsec:unitary}

A \emph{unitary fusion category} is a fusion category with a $C^*$-structure such that the associativity and unit constraints are unitary and $(f\otimes g)^*=f^*\otimes g^*$ for all morphisms $f$ and $g$. On objects the superscript $*$ denotes the dual, and on morphisms of a unitary fusion category it denotes the adjoint. A fusion category is \emph{unitarizable} if it is tensor equivalent to a unitary fusion category. Both properties pass to fusion subcategories. A \emph{unitary braided fusion category} is a unitary fusion category with a braiding. The braiding is then unitary \cite[Theorem~3.2]{Gal14}.

For a simple object $Y$ and an object $V$, the $C^*$-structure gives $\Hom_\cC(Y,V)$ the inner product defined by
\begin{equation}\label{eq:inner-product}
u^*v=\langle u,v\rangle\id_Y ,
\end{equation}
linear in $v$. For a Hilbert space $L$ we take the object $L\odot Y$ with $\iota_\ell^*\,\iota_{\ell'}=\langle\ell,\ell'\rangle\,\id_Y$ for all $\ell,\ell'\in L$. With the inner products \eqref{eq:inner-product}, the isomorphism \eqref{eq:canonical-decomposition} is unitary. We call it the orthogonal decomposition of $V$ into simple objects.
These conventions for $C^*$-tensor categories follow \cite[Sections~2--3]{LR97} and \cite[Section~2.4]{MUe03}.

\subsection{Localizations}\label{subsec:localization}

Let $W$ be a nonzero finite-dimensional vector space. Recall from the introduction that a Yang--Baxter operator on $W$ is an invertible $R\in\End(W\otimes W)$ such that $R_{12}R_{23}R_{12}=R_{23}R_{12}R_{23}$ on $W^{\otimes3}$, where $R_{12}:=R\otimes\id_W$ and $R_{23}:=\id_W\otimes R$. The corresponding representation is
\[
\rho^R_n\colon\C B_n\to\End(W^{\otimes n}),\qquad \rho^R_n(\sigma_i)=\id_W^{\otimes(i-1)}\otimes R\otimes\id_W^{\otimes(n-i-1)}.
\]

\begin{definition}\label{def:localization}
Let $X$ be an object of a braided fusion category. A \emph{localization} of $X$ is a pair $(W,R)$, where $W$ is a nonzero finite-dimensional vector space and $R$ is a Yang--Baxter operator on $W$, such that for every $n\geq2$ there is an injective algebra homomorphism $\phi_n\colon\rho^X_n(\C B_n)\to\End(W^{\otimes n})$ with $\phi_n\circ\rho^X_n=\rho^R_n$. Its \emph{dimension} is $\dim W$. A localization is \emph{unitary} if $W$ is a Hilbert space and $R$ is unitary.
\end{definition}

This is \cite[Definition~1.2]{GHR13}. For unitary braid group representations it recovers \cite[Definition~2.3 and Remark~2.5]{RW12}, and a unitary localization is a localization in the sense of \cite[Definition~5.2]{RW18}.

All the localizations in this paper are obtained from Lemma~\ref{lem:localization} below. It concerns families of unital algebra homomorphisms $\Phi_n\colon E_n\to\End(W^{\otimes n})$, $n\geq1$, that are \emph{compatible with tensor products}, that is, such that
\begin{equation}\label{eq:tc}
\Phi_{n+k}(f\otimes g)=\Phi_n(f)\otimes\Phi_k(g)
\end{equation}
for all $n,k\geq1$, $f\in E_n$ and $g\in E_k$. This is the property \eqref{eq:tensor-compatible} of Theorem~\ref{thm:main}.

\begin{lemma}\label{lem:localization}
Let $X$ be an object of a braided fusion category $\cC$, let $W$ be a nonzero finite-dimensional vector space, and let $\Phi_n\colon E_n\to\End(W^{\otimes n})$, $n\geq1$, be a family of unital algebra homomorphisms that is compatible with tensor products. Put $R:=\Phi_2(c_{X,X})$.
\begin{enumerate}[label=\textup{(\alph*)}]
\item The operator $R$ is a Yang--Baxter operator on $W$, and $\Phi_n\circ\rho^X_n=\rho^R_n$ for every $n\geq2$.
\item If every $\Phi_n$ is injective, then $(W,R)$ is a localization of $X$.
\item If $\cC$ is unitary, $W$ is a Hilbert space and $\Phi_2(f^*)=\Phi_2(f)^*$ for every $f$, then $R$ is unitary.
\end{enumerate}
\end{lemma}

\begin{proof}
For $n\geq2$ and $1\leq i\leq n-1$, the morphism $\rho^X_n(\sigma_i)$ is the tensor product of $\id_X^{\otimes(i-1)}$, $c_{X,X}$ and $\id_X^{\otimes(n-i-1)}$. Applying \eqref{eq:tc} to this tensor product, and using that each $\Phi_k$ is unital, we get
\begin{equation}\label{eq:generators}
\Phi_n\bigl(\rho^X_n(\sigma_i)\bigr)
=\id_W^{\otimes(i-1)}\otimes R\otimes\id_W^{\otimes(n-i-1)} .
\end{equation}
Since $\Phi_2$ and $\Phi_3\circ\rho^X_3$ are unital algebra homomorphisms, $R$ is invertible, and the relation $\sigma_1\sigma_2\sigma_1=\sigma_2\sigma_1\sigma_2$ in $B_3$ gives the Yang--Baxter equation. By \eqref{eq:generators}, the algebra homomorphisms $\Phi_n\circ\rho^X_n$ and $\rho^R_n$ agree on the generators of $B_n$, and hence they are equal. This proves (a). If $\Phi_n$ is injective, then its restriction to the image of $\rho^X_n$ is the homomorphism $\phi_n$ required in Definition~\ref{def:localization}, which proves (b). For (c), the braiding $c_{X,X}$ is unitary \cite[Theorem~3.2]{Gal14}, and a unital algebra homomorphism that preserves adjoints sends unitary elements to unitary operators.
\end{proof}

%======================================================================
\section{The regular object of an integral fusion category}\label{sec:regular}

In this section $\cE$ is an integral fusion category with regular object $R_\cE$, and
\[
H(V):=\Hom_\cE(R_\cE,V),\qquad H(a)(g):=ag .
\]
No braiding is needed. For a nonzero object $T$ of $\cE$ we realize every algebra $\End_\cE(T^{\otimes n})$ on $W^{\otimes n}$, where $W:=H(T\otimes R_\cE)$ (Corollary~\ref{cor:regular}). The point is to join a morphism $R_\cE\to A\otimes R_\cE^{\otimes n}$ to a morphism $R_\cE\to B\otimes R_\cE^{\otimes k}$ and still have $n+k$ factors $R_\cE$ in the target. For this we use one morphism $\delta\colon R_\cE\to R_\cE\otimes R_\cE$. All diagrams are read from bottom to top, with the tensor factors ordered from left to right.

\begin{lemma}\label{lem:beta}
The functor $H$ and the regular object $R_\cE$ have the following properties.
\begin{enumerate}[label=\textup{(\alph*)}]
\item The functor $H\colon\cE\to\Vect$ is faithful, and $\dim H(V)=d_V$ for every object $V$.
\item There is a morphism $\delta\colon R_\cE\to R_\cE\otimes R_\cE$ such that, for all objects $S$ and $V$ of $\cE$, the following map is an isomorphism.
\begin{equation}\label{eq:Gamma}
\Gamma_{S,V}\colon H(V)\otimes\Hom_\cE(S,R_\cE)\longrightarrow\Hom_\cE(S,V\otimes R_\cE),\qquad g\otimes s\longmapsto(g\otimes\id_{R_\cE})\,\delta\,s .
\end{equation}
\item Suppose that $\cE$ is unitary, and let $\omega$ be a faithful positive linear functional on $\End_\cE(R_\cE)$. Then $\langle g,g'\rangle:=\omega(g^*g')$ is an inner product on $H(V)$, and $H(a)^*=H(a^*)$ for every morphism $a$. The morphism $\delta$ in \textup{(b)} can be chosen such that, for every object $V$ and all $g,g'\in H(V)$,
\begin{equation}\label{eq:delta-unitary}
\delta^*\,(g^*g'\otimes\id_{R_\cE})\,\delta=\langle g,g'\rangle\,\id_{R_\cE} .
\end{equation}
\end{enumerate}
\end{lemma}

\begin{proof}
In this proof $i$ and $j$ denote simple objects in $\cO(\cE)$.

(a) The multiplicity of $i$ in $R_\cE$ is $d_i$. Thus $\dim H(V)=\sum_id_i[V:i]=d_V$. In particular $H(i)\neq0$ for every $i$, and hence $H$ is faithful \cite[Lemma~3.1]{MT08}.

(b) We first construct an isomorphism between the two sides of \eqref{eq:Gamma}, regarded as functors of the pair $(S,V)$. They are linear in each variable, contravariant in $S$ and covariant in $V$. On a pair of simple objects $S=j$ and $V=i$ both have dimension $d_id_j$, since $i\otimes R_\cE$ is the direct sum of $d_i$ copies of $R_\cE$ by \eqref{eq:regular-absorbs}. Choose a linear isomorphism for each such pair. These choices extend uniquely to a natural isomorphism $\Gamma$ by semisimplicity. Indeed, the decomposition \eqref{eq:canonical-decomposition} is natural in $V$, and the same decomposition of $S$ gives the dual multiplicity spaces for a contravariant functor. The same argument is used in the proof of \cite[Proposition~6.1.14]{EGNO15}.

Now evaluate $\Gamma_{R_\cE,R_\cE}$ at $\id_{R_\cE}\otimes\id_{R_\cE}$, and call its value $\delta$. Naturality under $g\colon R_\cE\to V$ in the second variable and under $s\colon S\to R_\cE$ in the first gives exactly \eqref{eq:Gamma} (Figure~\ref{fig:gamma}). In particular, the single morphism $\delta$ determines the whole family.

\begin{figure}[H]
\centering
\begin{tikzpicture}[x=.9cm,y=.6cm]
  \draw[gwire] (1.15,-1.75) -- (1.15,1.00);
  \draw[gwire] (.98,1.00) .. controls (.98,1.45) and (.42,1.33) .. (.42,1.82) -- (.42,3.68);
  \draw[gwire] (1.32,1.00) .. controls (1.32,1.46) and (1.9,1.35) .. (1.9,1.86) -- (1.9,3.68);
  \node[gdelta] at (1.15,1.00) {$\delta$};
  \node[gbox] at (.42,2.45) {$g$};
  \node[gbox] at (1.15,-.62) {$s$};
  \node[glab,below] at (1.15,-1.75) {$S$};
  \node[glab,anchor=west] at (1.2,.2) {$R_\cE$};
  \node[glab,above] at (.42,3.68) {$V$};
  \node[glab,above] at (1.9,3.68) {$R_\cE$};
\end{tikzpicture}
\caption{The morphism $\Gamma_{S,V}(g\otimes s)=(g\otimes\id_{R_\cE})\,\delta\,s$ of \eqref{eq:Gamma}, read from bottom to top. The morphisms $g$ and $s$ are attached to $\delta$.}
\label{fig:gamma}
\end{figure}
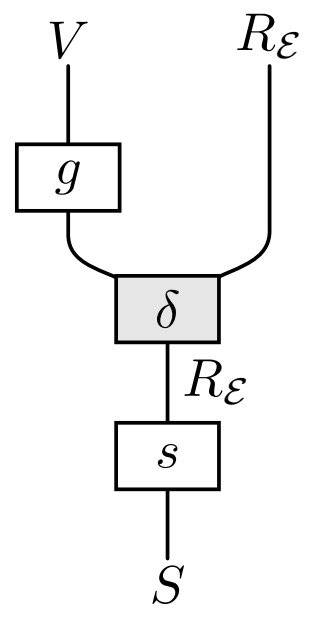

(c) Positive definiteness follows from the faithfulness of $\omega$, and $\langle ag,g'\rangle=\omega(g^*a^*g')=\langle g,a^*g'\rangle$ gives $H(a)^*=H(a^*)$. Give the spaces $\Hom_\cE(j,A)$ the inner products \eqref{eq:inner-product}, defined by $s^*t=\langle s,t\rangle\,\id_j$, and choose the isomorphisms $\Gamma_{j,i}$ unitary. For fixed $j$, both sides of \eqref{eq:Gamma} are functors of $V$ that preserve adjoints. Hence the orthogonal decomposition of $V$ into simple objects makes $\Gamma_{j,V}$ unitary as well. Testing this isometry on $g\otimes s$ and $g'\otimes t$ gives
\[
s^*\,\delta^*(g^*g'\otimes\id_{R_\cE})\,\delta\,t=\langle g,g'\rangle\,s^*t\qquad(s,t\colon j\to R_\cE).
\]
These matrix coefficients, for all $j$, determine an endomorphism of $R_\cE$. This proves \eqref{eq:delta-unitary}.
\end{proof}

The property in Lemma~\ref{lem:beta}(b) characterizes $R_\cE$ up to isomorphism. If a nonzero object $P$ has this property, then for every simple object $i$ the object $i\otimes P$ is the direct sum of $\dim\Hom_\cE(P,i)$ copies of $P$. Taking dimensions gives $d_id_P=[P:i]\,d_P$, and hence $P\cong R_\cE$.

We fix $\delta$ as in Lemma~\ref{lem:beta}(b). When $\cE$ is unitary we first choose a faithful positive linear functional $\omega$ on the finite-dimensional $C^*$-algebra $\End_\cE(R_\cE)$, and we take $\delta$ as in Lemma~\ref{lem:beta}(c). For $g\colon R_\cE\to V$ put
\begin{equation}\label{eq:hat}
\widehat g:=(g\otimes\id_{R_\cE})\,\delta\colon R_\cE\longrightarrow V\otimes R_\cE .
\end{equation}
In particular, $\widehat{\id_{R_\cE}}=\delta$. For $u\colon V\to V'$ the definition gives
\begin{equation}\label{eq:hat-natural}
\widehat{ug}=(u\otimes\id_{R_\cE})\,\widehat g .
\end{equation}
This is Figure~\ref{fig:hat}(a) with $u$ attached above the strand $V$. Thus a hat is simply the factor attached on the left of $s$ in \eqref{eq:Gamma}. Equivalently, as $S$ varies, that isomorphism expresses $V\otimes R_\cE\cong H(V)\odot R_\cE$ by the Yoneda lemma. The copy map $g\odot\id_{R_\cE}$ followed by this identification is $\widehat g$.

\medskip\noindent\emph{Joining two blocks.} Write $\id_m:=\id_{R_\cE^{\otimes m}}$, with $R_\cE^{\otimes0}=\one$. For $n\geq1$ set $F_n(A):=H(A\otimes R_\cE^{\otimes n})$ and $F_n(a):=H(a\otimes\id_n)$. Given $f\in F_n(A)$ and $g\in F_k(B)$, attach $\widehat g$ to the first strand $R_\cE$ that leaves $f$. Its $k+1$ strands $R_\cE$ replace that one strand. This defines a linear map
\begin{equation}\label{eq:J}
J_{n,k}\colon F_n(A)\otimes F_k(B)\longrightarrow F_{n+k}(A\otimes B),\qquad f\otimes g\longmapsto f\star g:=\bigl(\id_A\otimes\widehat g\otimes\id_{n-1}\bigr)f .
\end{equation}
Figure~\ref{fig:hat}(b) shows the morphism $f\star g$. The first strand $R_\cE$ that is available for a further insertion is inside the block $R_\cE^{\otimes k}$ that comes from $g$.

\begin{figure}[H]
\centering
\begin{tikzpicture}[x=.93cm,y=.77cm]
  \draw[gwire] (1.15,0) -- (1.15,1.00);
  \draw[gwire] (.98,1.00) .. controls (.98,1.45) and (.42,1.33) .. (.42,1.82) -- (.42,3.68);
  \draw[gwire] (1.32,1.00) .. controls (1.32,1.46) and (1.9,1.35) .. (1.9,1.86) -- (1.9,3.68);
  \node[gdelta] at (1.15,1.00) {$\delta$};
  \node[gbox] at (.42,2.45) {$g$};
  \node[glab,below] at (1.15,0) {$R_\cE$};
  \node[glab,above] at (.42,3.68) {$V$};
  \node[glab,above] at (1.9,3.68) {$R_\cE$};
  \node[glab,anchor=north] at (1.15,-.62) {(a) $\widehat g$};
\end{tikzpicture}%
\hspace{1.5cm}%
\begin{tikzpicture}[x=1.3cm,y=.77cm]
  \draw[gwire] (2.45,0) -- (2.45,.64);
  \draw[gwire] (.42,.64) -- (.42,3.68);
  \draw[gwire] (2.35,.64) -- (2.35,1.62);
  \draw[gwire] (4.52,.64) -- (4.52,3.68);
  \draw[gwire] (2.22,1.62) .. controls (2.22,2.03) and (1.96,2.0) .. (1.96,2.46);
  \draw[gwire] (2.48,1.62) .. controls (2.48,2.07) and (3.35,1.98) .. (3.35,2.50) -- (3.35,3.68);
  \draw[gwire] (1.49,2.46) -- (1.49,3.68);
  \draw[gwire] (2.46,2.46) -- (2.46,3.68);
  \node[gbox,minimum width=5.7cm] at (2.47,.64) {$f$};
  \node[gdelta] at (2.35,1.62) {$\delta$};
  \node[gbox,minimum width=1.75cm] at (1.975,2.46) {$g$};
  \node[glab,below] at (2.45,0) {$R_\cE$};
  \foreach \x/\s in {.42/{A},1.49/{B},2.46/{R_\cE^{\otimes k}},3.35/{R_\cE},4.52/{R_\cE^{\otimes(n-1)}}}
      \node[glab,above] at (\x,3.68) {$\s$};
  \node[glab,anchor=north] at (2.45,-.62) {(b) $f\star g$};
\end{tikzpicture}
\caption{(a) The morphism $\widehat g=(g\otimes\id_{R_\cE})\,\delta$. (b) The morphism $f\star g$ for $f\in F_n(A)$ and $g\in F_k(B)$. A strand labeled by a tensor power represents the whole block, and the last strand is omitted when $n=1$.}
\label{fig:hat}
\end{figure}
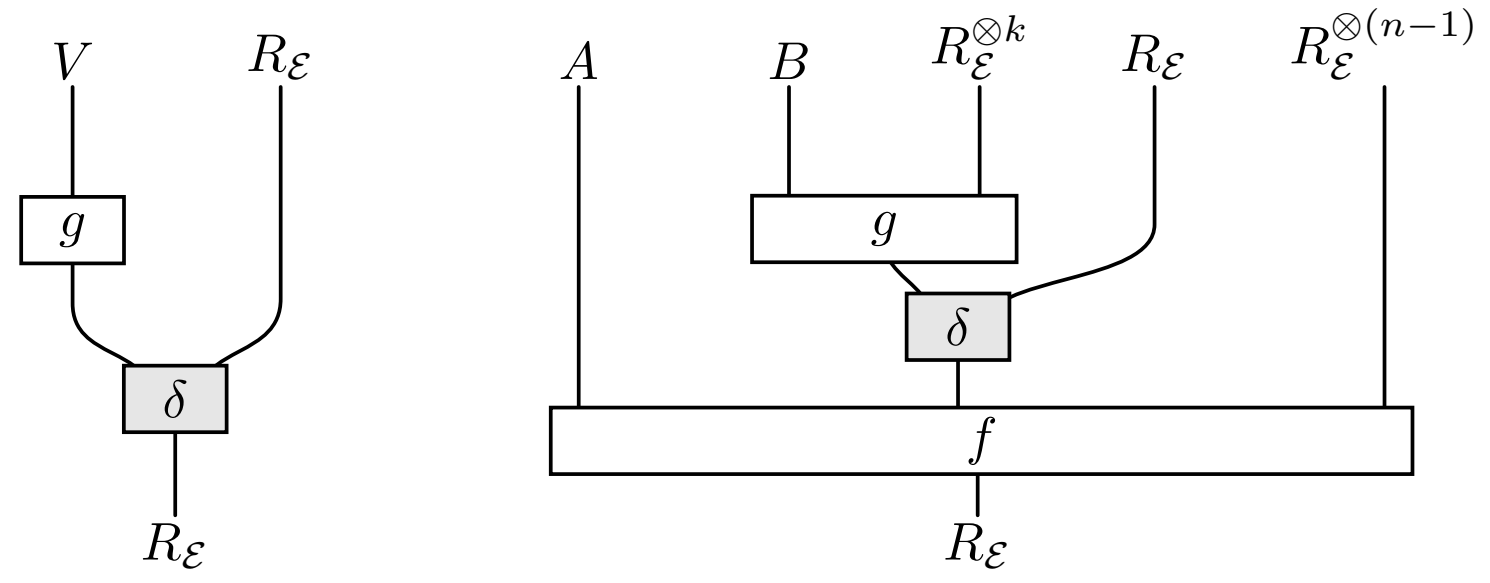
\begin{theorem}\label{thm:regular}
For $n,k,l\geq1$ and objects $A$, $B$, $C$ of $\cE$, the maps $J_{n,k}$ in \eqref{eq:J} satisfy the following.
\begin{enumerate}[label=\textup{(\alph*)}]
\item Each $J_{n,k}$ is an isomorphism of vector spaces.
\item They are natural in $A$ and $B$. For $a\colon A\to A'$ and $b\colon B\to B'$,
\[
F_{n+k}(a\otimes b)\,J_{n,k}=J_{n,k}\,\bigl(F_n(a)\otimes F_k(b)\bigr).
\]
\item They are associative. For $f\in F_n(A)$, $g\in F_k(B)$ and $h\in F_l(C)$,
\[
(f\star g)\star h=f\star(g\star h).
\]
\item If $\cE$ is unitary, then, with the choices of Lemma~\textup{\ref{lem:beta}(c)}, every $J_{n,k}$ is unitary and $F_n(a^*)=F_n(a)^*$.
\end{enumerate}
\end{theorem}

\begin{proof}
(a) The isomorphism in Lemma~\ref{lem:beta}(b) already gives $J_{n,k}$. To see it, put $D:=R_\cE^{\otimes(n-1)}$, $V:=B\otimes R_\cE^{\otimes k}$ and $S:=A^*\otimes R_\cE\otimes{}^*D$. Bending the two outer strands, as in Figure~\ref{fig:mate}, gives an isomorphism
\[
\Hom_\cE(R_\cE,A\otimes M\otimes D)\cong\Hom_\cE(S,M),
\]
natural in $M$ \cite[Proposition~2.10.8]{EGNO15}. We call the image of a morphism its mate. The caps in the figure are the evaluations $A^*\otimes A\to\one$ and $D\otimes{}^*D\to\one$.
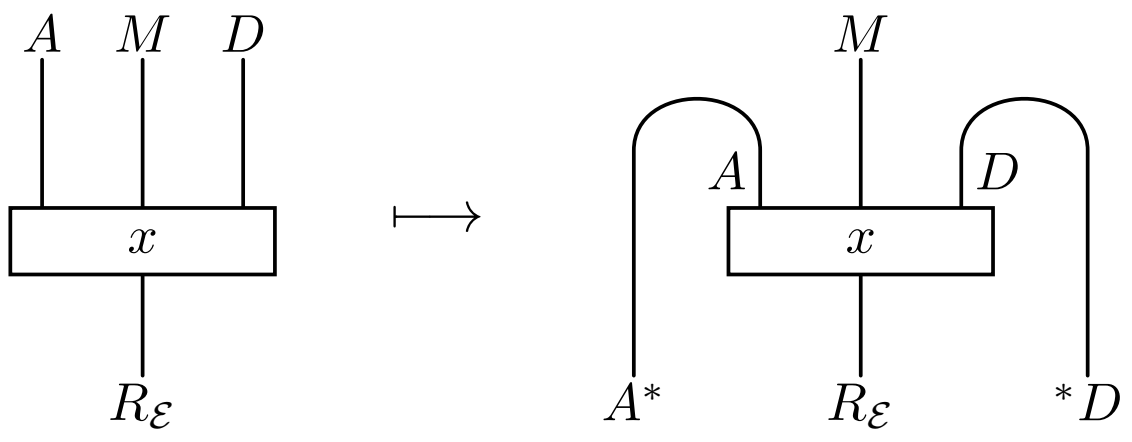
\begin{figure}[H]
\centering
\begin{tikzpicture}[x=.95cm,y=.76cm]
  \draw[gwire] (0,0) -- (0,1.17);
  \foreach \u/\lab in {-.70/A,0/M,.70/D}{
    \draw[gwire] (\u,1.17) -- (\u,2.75);
    \node[glab,above] at (\u,2.75) {$\lab$};}
  \node[gbox,minimum width=1.75cm] at (0,1.17) {$x$};
  \node[glab,below] at (0,0) {$R_\cE$};
  \node at (2.05,1.35) {$\longmapsto$};
  \begin{scope}[xshift=4.75cm]
    \draw[gwire] (0,0) -- (0,2.75);
    \draw[gwire] (-1.58,0) -- (-1.58,1.96)
        .. controls (-1.58,2.56) and (-.70,2.56) .. (-.70,1.96) -- (-.70,1.17);
    \draw[gwire] (.70,1.17) -- (.70,1.96)
        .. controls (.70,2.56) and (1.58,2.56) .. (1.58,1.96) -- (1.58,0);
    \node[gbox,minimum width=1.75cm] at (0,1.17) {$x$};
    \node[glab,above] at (0,2.75) {$M$};
    \node[glab,below] at (-1.58,0) {$A^*$};
    \node[glab,below] at (0,0) {$R_\cE$};
    \node[glab,below] at (1.58,0) {${}^*D$};
    \node[glab,anchor=east] at (-.74,1.78) {$A$};
    \node[glab,anchor=west] at (.74,1.78) {$D$};
  \end{scope}
\end{tikzpicture}
\caption{The mate of a morphism $x\colon R_\cE\to A\otimes M\otimes D$, obtained by bending the two outer strands.}
\label{fig:mate}
\end{figure}
Take $M=R_\cE$ and $M=V\otimes R_\cE$. Then $F_n(A)=\Hom_\cE(R_\cE,A\otimes R_\cE\otimes D)$, $F_k(B)=H(V)$ and $F_{n+k}(A\otimes B)=\Hom_\cE(R_\cE,A\otimes V\otimes R_\cE\otimes D)$. If $s\colon S\to R_\cE$ is the mate of $f$, then the mate of $f\star g=(\id_A\otimes\widehat g\otimes\id_D)f$ is $\widehat g\,s=\Gamma_{S,V}(g\otimes s)$, by naturality in $M$. Thus, after passing to mates and interchanging the two factors of the source, $J_{n,k}$ is the isomorphism $\Gamma_{S,V}$. Only the duals and \eqref{eq:Gamma} enter here.

(b) In Figure~\ref{fig:hat}(b), the morphism $a$ acts on the strand $A$ of $f$ and $b$ acts on the strand $B$ of $g$. Composing there with $a$ and $b$ gives the diagram of $F_n(a)(f)\star F_k(b)(g)$, by \eqref{eq:hat-natural} applied to $u=b\otimes\id_k$. This is the asserted identity.

(c) Let $f\circ_1g:=(\id_A\otimes g\otimes\id_{n-1})f$ be the morphism in which $g$ is attached directly to the first strand $R_\cE$ that leaves $f$. It is an element of $F_{n+k-1}(A\otimes B)$. This insertion is associative, and hence a product of several factors does not depend on the parentheses. Indeed, since $k\geq1$, the first strand $R_\cE$ of $f\circ_1g$ comes from $g$, and hence $(f\circ_1g)\circ_1h$ and $f\circ_1(g\circ_1h)$ have the same connections. Regarding $\delta$ as an element of $F_2(\one)$, we have
\[
f\star g=f\circ_1\delta\circ_1g .
\]
Both ways of placing the parentheses in the triple product therefore give $f\circ_1\delta\circ_1g\circ_1\delta\circ_1h$. Figure~\ref{fig:assoc} shows the case $n=k=l=1$.

\begin{figure}[H]
\centering
\begin{tikzpicture}[x=1.12cm,y=.72cm]
  \draw[gwire] (1.40,0) -- (1.40,.50);
  \draw[gwire] (.24,.50) -- (.24,4.48);
  \draw[gwire] (2.41,.50) -- (2.41,1.32);
  \draw[gwire] (2.28,1.32) .. controls (2.28,1.77) and (1.97,1.67) .. (1.97,2.14);
  \draw[gwire] (2.54,1.32) .. controls (2.54,1.83) and (4.82,1.54) .. (4.82,2.28) -- (4.82,4.48);
  \draw[gwire] (1.22,2.14) -- (1.22,4.48);
  \draw[gwire] (2.71,2.14) -- (2.71,2.94);
  \draw[gwire] (2.58,2.94) .. controls (2.58,3.36) and (2.49,3.31) .. (2.49,3.75);
  \draw[gwire] (2.84,2.94) .. controls (2.84,3.38) and (3.88,3.20) .. (3.88,3.85) -- (3.88,4.48);
  \draw[gwire] (2.07,3.75) -- (2.07,4.48);
  \draw[gwire] (2.94,3.75) -- (2.94,4.48);
  \node[gbox,minimum width=2.9cm] at (1.37,.50) {$f$};
  \node[gdelta] at (2.41,1.32) {$\delta$};
  \node[gbox,minimum width=2.05cm] at (1.97,2.14) {$g$};
  \node[gdelta] at (2.71,2.94) {$\delta$};
  \node[gbox,minimum width=1.35cm] at (2.49,3.75) {$h$};
  \node[glab,below] at (1.40,0) {$R_\cE$};
  \foreach \x/\s in {.24/A,1.22/B,2.07/C,2.94/R_\cE,3.88/R_\cE,4.82/R_\cE}
      \node[glab,above] at (\x,4.48) {$\s$};
\end{tikzpicture}
\caption{The morphism $(f\star g)\star h=f\star(g\star h)$ for $n=k=l=1$. The morphism $h$ is attached to the first strand $R_\cE$ that leaves $g$, and the box $g$ separates the two boxes $\delta$.}
\label{fig:assoc}
\end{figure}
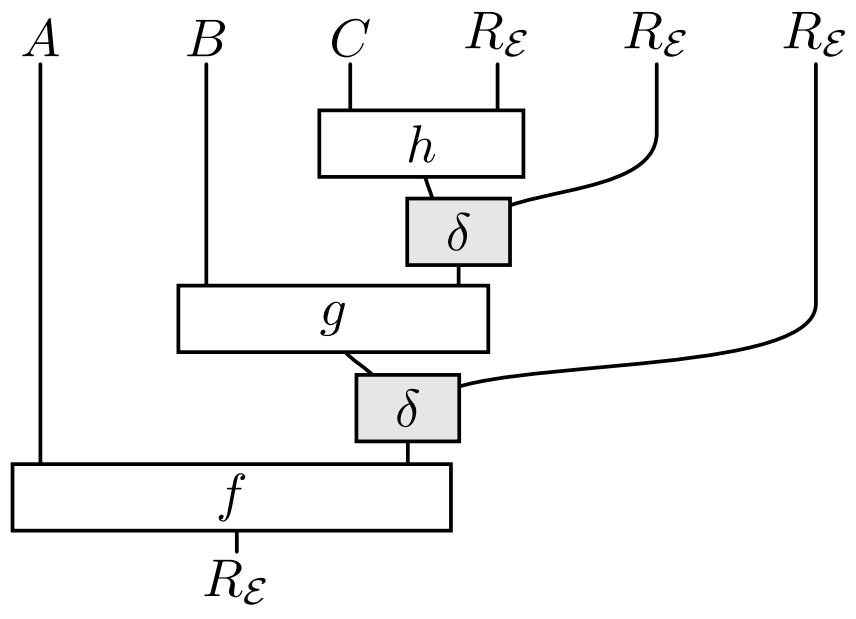

(d) By \eqref{eq:delta-unitary}, $\widehat g^{\,*}\,\widehat{g'}=\langle g,g'\rangle\,\id_{R_\cE}$ for $g,g'\in F_k(B)$. In Figure~\ref{fig:contract} the two boxes $\delta$ and $\delta^*$ contract around $a=g^*g'$, and $\omega(a)=\langle g,g'\rangle$.
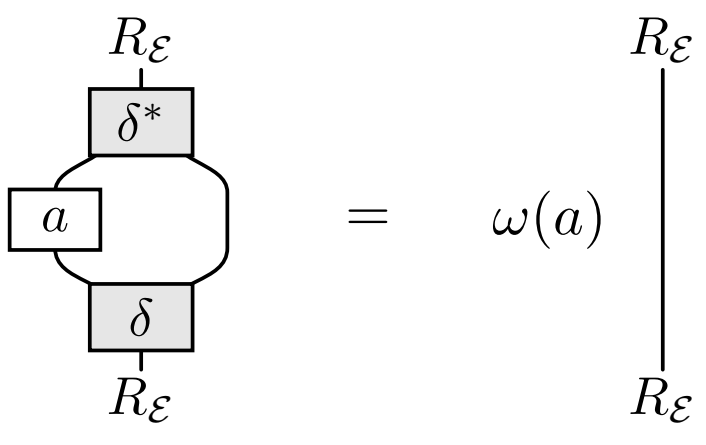
\begin{figure}[H]
\centering
\begin{tikzpicture}[x=1cm,y=.62cm]
  \draw[gwire] (0,0) -- (0,.56);
  \draw[gwire] (-.16,.56) .. controls (-.16,.98) and (-.57,.91) .. (-.57,1.30)
    -- (-.57,1.89) .. controls (-.57,2.22) and (-.16,2.20) .. (-.16,2.64);
  \draw[gwire] (.16,.56) .. controls (.16,.98) and (.57,.91) .. (.57,1.30)
    -- (.57,1.89) .. controls (.57,2.22) and (.16,2.20) .. (.16,2.64);
  \draw[gwire] (0,2.64) -- (0,3.20);
  \node[gdelta] at (0,.56) {$\delta$};
  \node[gdelta] at (0,2.64) {$\delta^*$};
  \node[gbox,minimum width=.6cm,minimum height=.4cm] at (-.57,1.60) {$a$};
  \node at (1.5,1.6) {$=$};
  \node at (2.7,1.6) {$\omega(a)$};
  \draw[gwire] (3.45,0) -- (3.45,3.20);
  \foreach \x in {0,3.45} {
    \node[glab,below] at (\x,0) {$R_\cE$};
    \node[glab,above] at (\x,3.20) {$R_\cE$};}
\end{tikzpicture}
\caption{The identity \eqref{eq:delta-unitary} for $a=g^*g'$, where $\omega(a)=\langle g,g'\rangle$.}
\label{fig:contract}
\end{figure}
Applying this contraction inside \eqref{eq:J} gives
\[
\langle f\star g,f'\star g'\rangle=\omega\bigl(f^*(\id_A\otimes\widehat g^{\,*}\,\widehat{g'}\otimes\id_{n-1})f'\bigr)=\langle f,f'\rangle\,\langle g,g'\rangle .
\]
Thus $J_{n,k}$ preserves inner products, and it is surjective by (a). Finally, $F_n(a)^*=H\bigl((a\otimes\id_n)^*\bigr)=F_n(a^*)$ by Lemma~\ref{lem:beta}(c).
\end{proof}

\medskip\noindent\emph{Tensor powers of a fixed object.} Fix a nonzero object $T$ of $\cE$, set $W:=F_1(T)=\Hom_\cE(R_\cE,T\otimes R_\cE)$, and define $J_n\colon W^{\otimes n}\to F_n(T^{\otimes n})$ recursively by
\[
J_1:=\id_W,\qquad J_{n+1}:=J_{n,1}\,(J_n\otimes\id_W).
\]
These maps are isomorphisms by Theorem~\ref{thm:regular}(a).

\begin{corollary}\label{cor:regular}
For $n\geq1$ and $f\in\End_\cE(T^{\otimes n})$ set $\Phi^T_n(f):=J_n^{-1}\,F_n(f)\,J_n$.
\begin{enumerate}[label=\textup{(\alph*)}]
\item The maps $\Phi^T_n\colon\End_\cE(T^{\otimes n})\to\End(W^{\otimes n})$ are injective unital algebra homomorphisms, and $\dim W=d_T\FPdim(\cE)$.
\item For $n,k\geq1$, $f\in\End_\cE(T^{\otimes n})$ and $g\in\End_\cE(T^{\otimes k})$, we have $\Phi^T_{n+k}(f\otimes g)=\Phi^T_n(f)\otimes\Phi^T_k(g)$.
\item If $\cE$ is unitary, then, with the choices of Lemma~\textup{\ref{lem:beta}(c)}, $\Phi^T_n(f^*)=\Phi^T_n(f)^*$.
\end{enumerate}
\end{corollary}

\begin{proof}
(a) The functor $-\otimes R_\cE^{\otimes n}$ is faithful by \eqref{eq:faithful-tensor}, and $H$ is faithful by Lemma~\ref{lem:beta}(a). Hence $F_n$ is faithful. Conjugating its action on endomorphisms by $J_n$ gives an injective unital algebra homomorphism. Lemma~\ref{lem:beta}(a) gives $\dim W=d_{T\otimes R_\cE}=d_T\sum_id_i^2=d_T\FPdim(\cE)$.

(b) On pure tensors, $J_n$ forms the iterated product $J_n(w_1\otimes\dots\otimes w_n)=(\cdots(w_1\star w_2)\star\cdots)\star w_n$. Associativity lets us group the first $n$ and the last $k$ factors, and hence
\begin{equation}\label{eq:iterated-J}
J_{n+k}=J_{n,k}\,(J_n\otimes J_k).
\end{equation}
Formally this is an induction on $k$, in which the case $k=1$ is the definition of $J_{n+1}$ and the inductive step is Theorem~\ref{thm:regular}(c) with $l=1$. With this factorization, the naturality of $J_{n,k}$ gives
\begin{align*}
F_{n+k}(f\otimes g)\,J_{n+k}&=J_{n,k}\,\bigl(F_n(f)\,J_n\otimes F_k(g)\,J_k\bigr)\\
&=J_{n+k}\,\bigl(\Phi^T_n(f)\otimes\Phi^T_k(g)\bigr).
\end{align*}
Composing with $J_{n+k}^{-1}$ on the left proves (b).

(c) Each $J_n$ is a composite of unitary maps, and $F_n$ preserves adjoints by Theorem~\ref{thm:regular}(d). The same is therefore true of $\Phi^T_n$.
\end{proof}

\begin{remark}\label{rem:regular-literature}
Suppose that $\delta$ is coassociative, that is, $(\delta\otimes\id_{R_\cE})\,\delta=(\id_{R_\cE}\otimes\delta)\,\delta$, and that it has a counit. Then $R_\cE$ is a right comodule over itself, and $\widehat g$ is the lift of $g$ to a morphism of comodules from $R_\cE$ to the cofree right $R_\cE$-comodule $V\otimes R_\cE$. The fact that $\widehat g$ is a morphism of comodules is precisely the step that uses coassociativity. The maps $J_{n,k}$ only use insertion into the first strand $R_\cE$, and for this reason they do not need it. In fact, the morphism $\delta$ of Lemma~\ref{lem:beta}(b) can be chosen coassociative if and only if $\cE$ admits a fiber functor. Indeed, suppose that $\delta$ is coassociative. Then the maps $H(V)\otimes H(V')\to H(V\otimes V')$, $g\otimes g'\mapsto(g\otimes g')\,\delta$, are natural and associative, and for $h:=(g\otimes g')\,\delta$ coassociativity gives $\widehat h=(\id_V\otimes\widehat{g'})\,\widehat g$. The map $g\mapsto\widehat g$ is injective by Lemma~\ref{lem:beta}(b). The map $x\otimes g'\mapsto(\id_V\otimes\widehat{g'})\,x$, for $x\in\Hom_\cE(R_\cE,V\otimes R_\cE)$, is also injective, because under the isomorphisms $\Hom_\cE(R_\cE,V\otimes M)\cong\Hom_\cE(V^*\otimes R_\cE,M)$ it becomes the map $\Gamma_{V^*\otimes R_\cE,V'}$ of \eqref{eq:Gamma}. Hence the maps $g\otimes g'\mapsto(g\otimes g')\,\delta$ are injective, and they are isomorphisms because both spaces have dimension $d_Vd_{V'}$. Moreover, $H(\one)$ has dimension $1$, and $H$ is exact and faithful. These maps are therefore a tensor structure on $H$ \cite[Definition~2.4.1]{EGNO15}, and $H$ is a fiber functor. The unit $\varepsilon$ of the one-dimensional algebra $H(\one)$ acts as the identity on every space $H(V)$, on both sides, and for $V=R_\cE$ this gives $(\varepsilon\otimes\id_{R_\cE})\,\delta=\id_{R_\cE}=(\id_{R_\cE}\otimes\varepsilon)\,\delta$. Thus a coassociative morphism $\delta$ with the property of Lemma~\ref{lem:beta}(b) always has a counit. Conversely, a fusion category with a fiber functor is tensor equivalent to the category of finite-dimensional modules over a semisimple Hopf algebra $K$ \cite[Theorem~5.3.12]{EGNO15}. Its regular object is $K$ with the action by left multiplication, and the comultiplication $\Delta\colon K\to K\otimes K$ is a coassociative morphism of $K$-modules with the property of Lemma~\ref{lem:beta}(b). Indeed, for a $K$-module $V$ let $V_0$ be the vector space $V$ with the trivial action. The map $v\otimes a\mapsto\Delta(a)(v\otimes1)$ is an isomorphism of $K$-modules from $V_0\otimes K$ onto $V\otimes K$, and composition with it is the map $\Gamma_{S,V}$ for every $S$, after $H(V)$ is identified with $V_0$ by $g\mapsto g(1)$.
\end{remark}

%======================================================================
\section{Dual strands and the proof of Theorem~\ref{thm:main}}\label{sec:sufficiency}

Throughout this section $X$ is a nonzero object of a braided fusion category $\cC$, and we put
\[
d:=d_X,\qquad \cD:=\langle X\rangle,\qquad \cA:=\langle X\otimes X^*\rangle .
\]
Recall that $E_n=\End_\cC(X^{\otimes n})$.
If $d^2$ is an integer, then $\cA$ is integral (Lemma~\ref{lem:integrality}), and Section~\ref{sec:regular} applies to the object $X\otimes X^*$ of $\cA$. To use it for $X$, we attach a copy of $X^*$ to each copy of $X$. Subsection~\ref{subsec:dual} constructs homomorphisms $S_n\colon E_n\to\End_\cA((X\otimes X^*)^{\otimes n})$ that are compatible with tensor products, and Subsection~\ref{subsec:sufficiency} composes them with the homomorphisms of Corollary~\ref{cor:regular}.

\subsection{The fusion subcategories generated by \texorpdfstring{$X$}{X} and by \texorpdfstring{$X\otimes X^*$}{X tensor X*}}\label{subsec:aux}

We first relate the fusion subcategory $\cD$ generated by $X$ to the fusion subcategory $\cA$ generated by $X\otimes X^*$. Besides the integrality of $\cA$, Lemma~\ref{lem:integrality} shows that $\cD$ has a canonical spherical structure when $d^2$ is an integer. This structure is used in the trace identity.

\begin{lemma}\label{lem:integrality}
Let $X$ be a nonzero object of a braided fusion category, with $d$, $\cD$ and $\cA$ as above.
\begin{enumerate}[label=\textup{(\alph*)}]
\item For every simple object $Y$ of $\cD$, the object $Y\otimes Y^*$ belongs to $\cA$.
\item Suppose that $d^2\in\Z$. Then $\cA$ is integral, $d_Y^2\in\Z$ for every simple object $Y$ of $\cD$, and $\cD$ is weakly integral. In particular, $\cD$ is pseudo-unitary, and it has a canonical spherical structure, with $\dim(Y)=d_Y$ for every simple object $Y$ of $\cD$.
\item If $X$ is simple, then $d^2\in\Z$ if and only if $\cD$ is weakly integral.
\end{enumerate}
\end{lemma}

\begin{proof}
(a) By Lemma~\ref{lem:generated}, $Y$ is a direct summand of $X^{\otimes r}$ for some $r\geq1$. Then $Y\otimes Y^*$ is a direct summand of $X^{\otimes r}\otimes(X^{\otimes r})^*$, and $(X^{\otimes r})^*\cong(X^*)^{\otimes r}$, since the dual of a tensor product is the tensor product of the duals in reverse order. By means of the braiding we can reorder the factors, and we get $X^{\otimes r}\otimes(X^*)^{\otimes r}\cong(X\otimes X^*)^{\otimes r}$. This object belongs to $\cA$, and therefore its direct summand $Y\otimes Y^*$ belongs to $\cA$.

(b) Let $i$ be a simple object of $\cA$. By Lemma~\ref{lem:generated}, applied to the object $X\otimes X^*$, it is a direct summand of $(X\otimes X^*)^{\otimes r}$ for some $r\geq1$, and the Frobenius--Perron dimension of this object is the integer $d^{2r}$. If $d_A+d_B$ is an integer for two objects $A$ and $B$, then $d_A$ and $d_B$ are integers \cite[Lemma~3.5.6]{EGNO15}. Hence $d_i$ is an integer, and $\cA$ is integral. By~\textup{(a)}, $d_Y^2=d_{Y\otimes Y^*}$ is an integer for every simple object $Y$ of $\cD$, and therefore $\FPdim(\cD)=\sum_{Y\in\cO(\cD)}d_Y^2$ is an integer. The assertions about the canonical spherical structure were recalled in Subsection~\ref{subsec:spherical}.

(c) One implication is contained in (b). Conversely, the square of the Frobenius--Perron dimension of every simple object of a weakly integral fusion category is an integer \cite[Proposition~9.6.9(i)]{EGNO15}.
\end{proof}

\subsection{Dual strands}\label{subsec:dual}
\begin{figure}[!ht]
\centering
\begin{tikzpicture}[x=1.0cm,y=1.6cm]
  % X strands
  \draw[thick] (0,0) -- (0,2);
  \draw[thick] (2,0) .. controls (2,1) and (1,1) .. (1,2);
  \draw[thick] (4,0) .. controls (4,1) and (2,1) .. (2,2);
  % dual strands over the X strands
  \draw[white,line width=5pt] (1,0) .. controls (1,1) and (3,1) .. (3,2);
  \draw[thick,dashed] (1,0) .. controls (1,1) and (3,1) .. (3,2);
  \draw[white,line width=5pt] (3,0) .. controls (3,1) and (4,1) .. (4,2);
  \draw[thick,dashed] (3,0) .. controls (3,1) and (4,1) .. (4,2);
  \draw[thick,dashed] (5,0) -- (5,2);
  \foreach \x/\l in {0/X,1/X^*,2/X,3/X^*,4/X,5/X^*}{\node[below] at (\x,0) {$\l$};}
  \foreach \x/\l in {0/X,1/X,2/X,3/X^*,4/X^*,5/X^*}{\node[above] at (\x,2) {$\l$};}
\end{tikzpicture}
\caption{The isomorphism $s_3\colon (X\otimes X^*)^{\otimes3}\to X^{\otimes3}\otimes(X^*)^{\otimes3}$, read from bottom to top. The dual strands are dashed, and each crossing is a braiding $c_{X^*,X}$. We draw $c_{V,V'}$ with the strand of $V$ passing over the strand of $V'$.}
\label{fig:dual}
\end{figure}
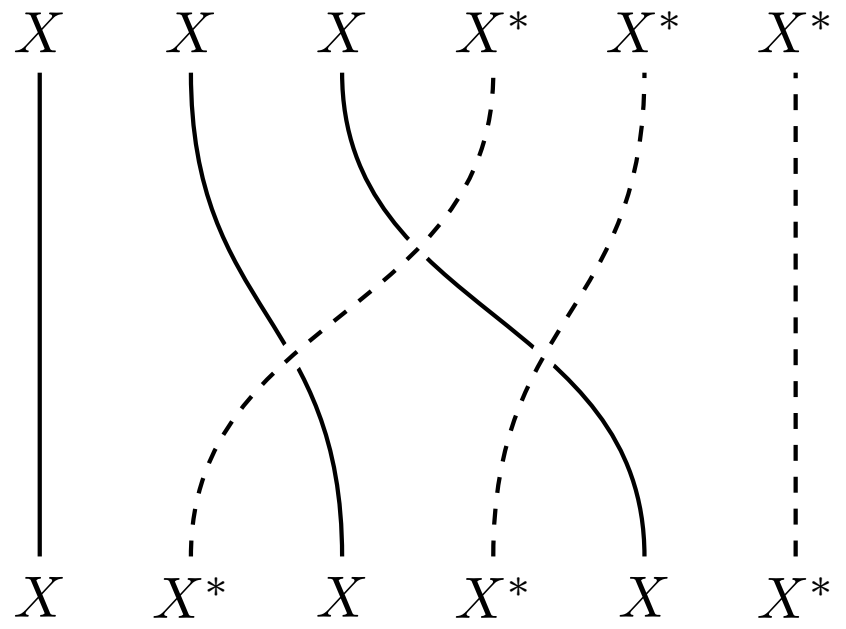

We define isomorphisms $s_n\colon (X\otimes X^*)^{\otimes n}\to X^{\otimes n}\otimes (X^*)^{\otimes n}$ by $s_1:=\id_{X\otimes X^*}$ and
\begin{equation}\label{eq:sn}
s_{n+1}:=(\id_{X^{\otimes n}}\otimes c_{(X^*)^{\otimes n},X}\otimes\id_{X^*})(s_n\otimes\id_{X\otimes X^*}).
\end{equation}
Thus $s_n$ moves each copy of $X^*$ to the right, across the copies of $X$ that follow it, and it preserves the order of the copies of $X$ and the order of the copies of $X^*$ (Figure~\ref{fig:dual}). Lemma~\ref{lem:dual-strands} describes how to join two groups of strands that have already been put in this order.

\begin{lemma}\label{lem:dual-strands}
For $n,k\geq 1$ put
\[
b_{n,k}:=\id_{X^{\otimes n}}\otimes c_{(X^*)^{\otimes n},X^{\otimes k}}\otimes\id_{(X^*)^{\otimes k}}
  \colon X^{\otimes n}\otimes (X^*)^{\otimes n}\otimes X^{\otimes k}\otimes (X^*)^{\otimes k}
       \longrightarrow X^{\otimes(n+k)}\otimes (X^*)^{\otimes(n+k)}.
\]
Then $s_{n+k}=b_{n,k}\,(s_n\otimes s_k)$.
\end{lemma}

\begin{proof}
Let $\cC\times\cC$ be the category of pairs of objects of $\cC$, with the componentwise tensor product. The functor $M\colon\cC\times\cC\to\cC$, $M(A,B):=A\otimes B$, is a monoidal functor with the structure isomorphisms
\[
\id_A\otimes c_{B,C}\otimes\id_D\colon M(A,B)\otimes M(C,D)\longrightarrow M(A\otimes C,B\otimes D).
\]
Indeed, these isomorphisms are natural, and for three pairs $(A,B)$, $(C,D)$ and $(E,F)$ the two sides of the associativity axiom contain the braidings $c_{B\otimes D,E}$ and $c_{B,C\otimes E}$. Expanding them by \eqref{eq:hexagons}, each side becomes the composite of $c_{B,E}$ with $c_{B,C}$ and $c_{D,E}$, and the last two morphisms commute because they act on disjoint tensor factors. By \eqref{eq:sn}, the iterated structure isomorphism of $M$ at the object $(X,X^*)$ is $s_n$, and its structure isomorphism at the pair of objects $(X^{\otimes n},(X^*)^{\otimes n})$ and $(X^{\otimes k},(X^*)^{\otimes k})$ is $b_{n,k}$. Hence the assertion is the identity \eqref{eq:gamma-iterated} for $M$.
\end{proof}

We use $s_n$ to let an endomorphism of $X^{\otimes n}$ act on the copies of $X$ in $(X\otimes X^*)^{\otimes n}$.

\begin{proposition}\label{prop:dual-strands}
For $n\geq1$ and $f\in E_n$ put
\[
S_n(f):=s_n^{-1}\,\bigl(f\otimes\id_{(X^*)^{\otimes n}}\bigr)\,s_n\in\End_\cA((X\otimes X^*)^{\otimes n}).
\]
\begin{enumerate}[label=\textup{(\alph*)}]
\item The maps $S_n\colon E_n\to\End_\cA((X\otimes X^*)^{\otimes n})$ are injective unital algebra homomorphisms.
\item For all $n,k\geq1$, $f\in E_n$ and $g\in E_k$, we have $S_{n+k}(f\otimes g)=S_n(f)\otimes S_k(g)$.
\item If $\cC$ is unitary, then $S_n(f^*)=S_n(f)^*$.
\end{enumerate}
\end{proposition}

\begin{proof}
(a) The map $f\mapsto f\otimes\id_{(X^*)^{\otimes n}}$ is a unital algebra homomorphism, and it is injective by \eqref{eq:faithful-tensor}. Conjugation by the isomorphism $s_n$ is an algebra isomorphism from $\End_\cC(X^{\otimes n}\otimes (X^*)^{\otimes n})$ onto $\End_\cC((X\otimes X^*)^{\otimes n})$. The latter algebra is $\End_\cA((X\otimes X^*)^{\otimes n})$ because $\cA$ is a full subcategory.

(b) By the naturality of the braiding $c_{(X^*)^{\otimes n},X^{\otimes k}}$ in its second variable, the morphism $g$ slides through the crossing of $b_{n,k}$, and $f$ does not meet it (Figure~\ref{fig:slide}).
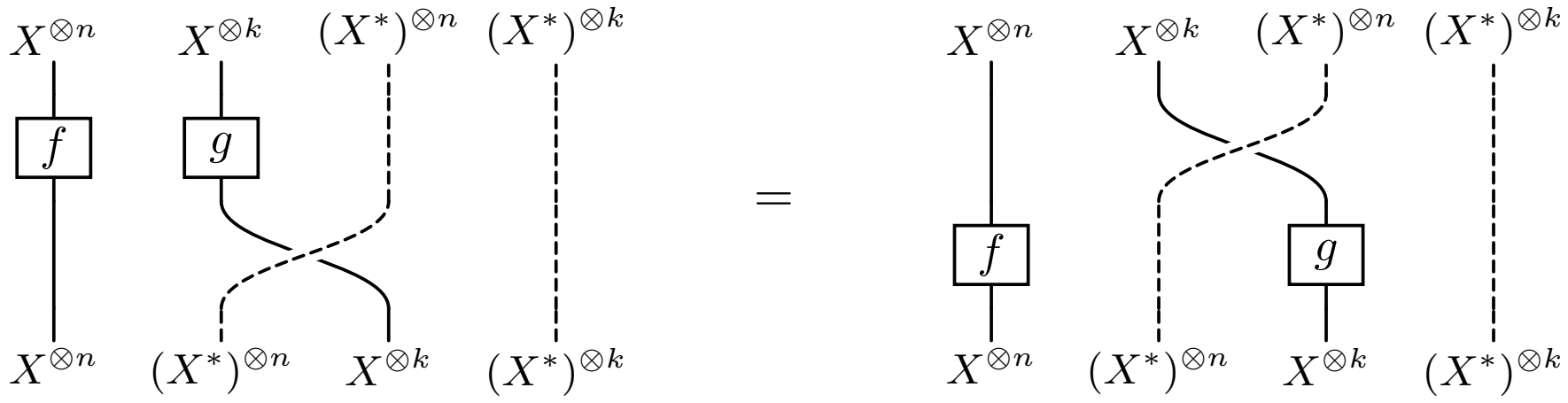
\begin{figure}[H]
\centering
\begin{tikzpicture}[x=1.25cm,y=.8cm,glab/.append style={font=\footnotesize}]
  \draw[gwire] (0.00,-.3)--(0.00,0);
  \node[glab,below] at (0.00,-.3) {$X^{\otimes n}$};
  \draw[gdash] (1.00,-.3)--(1.00,0);
  \node[glab,below] at (1.00,-.3) {$(X^*)^{\otimes n}$};
  \draw[gwire] (2.00,-.3)--(2.00,0);
  \node[glab,below] at (2.00,-.3) {$X^{\otimes k}$};
  \draw[gdash] (3.00,-.3)--(3.00,0);
  \node[glab,below] at (3.00,-.3) {$(X^*)^{\otimes k}$};
  \draw[gwire] (0.00,0)--(0.00,1);
  \draw[gdash] (3.00,0)--(3.00,1);
  \draw[gwire] (2.00,0).. controls (2.00,0.44) and (1.00,0.56)..(1.00,1);
  \draw[gdash,gover] (1.00,0).. controls (1.00,0.44) and (2.00,0.56)..(2.00,1);
  \draw[gwire] (0.00,1)--(0.00,2);
  \draw[gwire] (1.00,1)--(1.00,2);
  \draw[gdash] (2.00,1)--(2.00,2);
  \draw[gdash] (3.00,1)--(3.00,2);
  \draw[gwire] (0.00,2)--(0.00,2.3);
  \node[glab,above] at (0.00,2.3) {$X^{\otimes n}$};
  \draw[gwire] (1.00,2)--(1.00,2.3);
  \node[glab,above] at (1.00,2.3) {$X^{\otimes k}$};
  \draw[gdash] (2.00,2)--(2.00,2.3);
  \node[glab,above] at (2.00,2.3) {$(X^*)^{\otimes n}$};
  \draw[gdash] (3.00,2)--(3.00,2.3);
  \node[glab,above] at (3.00,2.3) {$(X^*)^{\otimes k}$};
  \node[gbox,minimum width=.55cm] at (0.00,1.5) {$f$};
  \node[gbox,minimum width=.55cm] at (1.00,1.5) {$g$};
  \node at (4.3,1) {$=$};
  \draw[gwire] (5.60,-.3)--(5.60,0);
  \node[glab,below] at (5.60,-.3) {$X^{\otimes n}$};
  \draw[gdash] (6.60,-.3)--(6.60,0);
  \node[glab,below] at (6.60,-.3) {$(X^*)^{\otimes n}$};
  \draw[gwire] (7.60,-.3)--(7.60,0);
  \node[glab,below] at (7.60,-.3) {$X^{\otimes k}$};
  \draw[gdash] (8.60,-.3)--(8.60,0);
  \node[glab,below] at (8.60,-.3) {$(X^*)^{\otimes k}$};
  \draw[gwire] (5.60,0)--(5.60,1);
  \draw[gdash] (6.60,0)--(6.60,1);
  \draw[gwire] (7.60,0)--(7.60,1);
  \draw[gdash] (8.60,0)--(8.60,1);
  \draw[gwire] (5.60,1)--(5.60,2);
  \draw[gdash] (8.60,1)--(8.60,2);
  \draw[gwire] (7.60,1).. controls (7.60,1.44) and (6.60,1.56)..(6.60,2);
  \draw[gdash,gover] (6.60,1).. controls (6.60,1.44) and (7.60,1.56)..(7.60,2);
  \draw[gwire] (5.60,2)--(5.60,2.3);
  \node[glab,above] at (5.60,2.3) {$X^{\otimes n}$};
  \draw[gwire] (6.60,2)--(6.60,2.3);
  \node[glab,above] at (6.60,2.3) {$X^{\otimes k}$};
  \draw[gdash] (7.60,2)--(7.60,2.3);
  \node[glab,above] at (7.60,2.3) {$(X^*)^{\otimes n}$};
  \draw[gdash] (8.60,2)--(8.60,2.3);
  \node[glab,above] at (8.60,2.3) {$(X^*)^{\otimes k}$};
  \node[gbox,minimum width=.55cm] at (5.60,.5) {$f$};
  \node[gbox,minimum width=.55cm] at (7.60,.5) {$g$};
\end{tikzpicture}
\caption{The naturality of the braiding $c_{(X^*)^{\otimes n},X^{\otimes k}}$. The morphism $g$ slides through the crossing of $b_{n,k}$.}
\label{fig:slide}
\end{figure}
In a formula,
\[
 (f\otimes g\otimes\id_{(X^*)^{\otimes(n+k)}})\,b_{n,k}
   = b_{n,k}\,\bigl(f\otimes\id_{(X^*)^{\otimes n}}\otimes g\otimes\id_{(X^*)^{\otimes k}}\bigr).
\]
By Lemma~\ref{lem:dual-strands} and this identity,
\begin{align*}
S_{n+k}(f\otimes g)&=(s_n\otimes s_k)^{-1}\,b_{n,k}^{-1}\,(f\otimes g\otimes\id_{(X^*)^{\otimes(n+k)}})\,b_{n,k}\,(s_n\otimes s_k)\\
&=(s_n\otimes s_k)^{-1}\bigl(f\otimes\id_{(X^*)^{\otimes n}}\otimes g\otimes\id_{(X^*)^{\otimes k}}\bigr)(s_n\otimes s_k),
\end{align*}
and the right-hand side is $S_n(f)\otimes S_k(g)$.

(c) By \eqref{eq:sn} and induction on $n$, the isomorphism $s_n$ is a composite of tensor products of braidings and identity morphisms. Every braiding of a unitary fusion category is unitary \cite[Theorem~3.2]{Gal14}, and tensor products and composites of unitary morphisms are unitary, because $(f\otimes g)^*=f^*\otimes g^*$. Hence $s_n$ is unitary, that is, $s_n^{-1}=s_n^*$. Therefore
\[
S_n(f)^*=s_n^*\,(f\otimes\id_{(X^*)^{\otimes n}})^*\,(s_n^{-1})^*=s_n^{-1}(f^*\otimes\id_{(X^*)^{\otimes n}})s_n=S_n(f^*).\qedhere
\]
\end{proof}

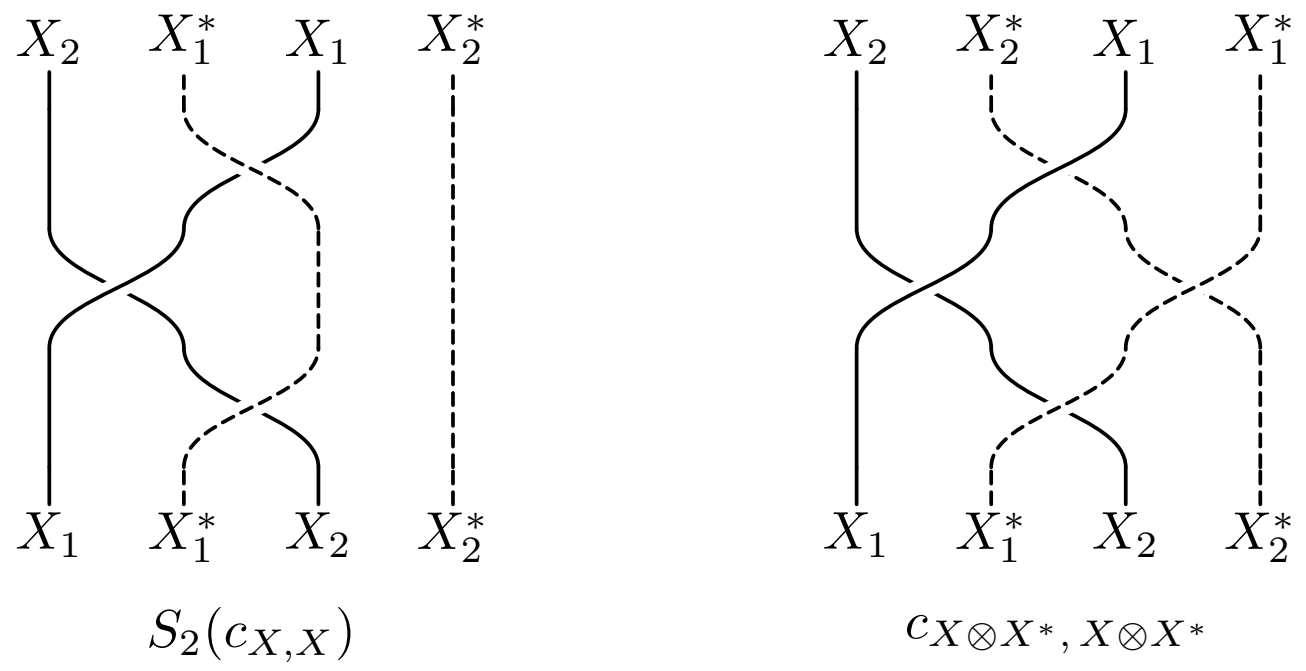
\begin{figure}[H]
\centering
\begin{tikzpicture}[x=.89cm,y=.79cm]
  \def\pgpos#1#2#3#4{%
    \begin{scope}[shift={(#1,#2)}]
      \draw[#4] (1,0).. controls (1,.44) and (0,.56)..(0,1);
      \draw[#3,gover] (0,0).. controls (0,.44) and (1,.56)..(1,1);
    \end{scope}}
  \def\pgneg#1#2#3#4{%
    \begin{scope}[shift={(#1,#2)}]
      \draw[#3] (0,0).. controls (0,.44) and (1,.56)..(1,1);
      \draw[#4,gover] (1,0).. controls (1,.44) and (0,.56)..(0,1);
    \end{scope}}
  \def\pgincoming{%
    \foreach \a/\sty/\name in {0/gwire/{X_1},1/gdash/{X_1^*},2/gwire/{X_2},3/gdash/{X_2^*}}{
      \draw[\sty] (\a,-.31)--(\a,0);
      \node[glab,below] at (\a,-.31) {$\name$};}}
  \def\pgoutgoing#1{%
    \foreach \a/\sty/\name in {#1}{
      \draw[\sty] (\a,3)--(\a,3.31);
      \node[glab,above] at (\a,3.31) {$\name$};}}
  \begin{scope}
    \pgincoming
    \draw[gwire] (0,0)--(0,1);
    \draw[gdash] (3,0)--(3,3);
    \pgpos{1}{0}{gdash}{gwire}
    \pgpos{0}{1}{gwire}{gwire}
    \draw[gdash] (2,1)--(2,2);
    \draw[gwire] (0,2)--(0,3);
    \pgneg{1}{2}{gwire}{gdash}
    \pgoutgoing{0/gwire/{X_2},1/gdash/{X_1^*},2/gwire/{X_1},3/gdash/{X_2^*}}
    \node[glab] at (1.5,-1.4) {$S_2(c_{X,X})$};
  \end{scope}
  \begin{scope}[shift={(6,0)}]
    \pgincoming
    \draw[gwire] (0,0)--(0,1);
    \draw[gdash] (3,0)--(3,1);
    \pgpos{1}{0}{gdash}{gwire}
    \pgpos{0}{1}{gwire}{gwire}
    \pgpos{2}{1}{gdash}{gdash}
    \draw[gwire] (0,2)--(0,3);
    \draw[gdash] (3,2)--(3,3);
    \pgpos{1}{2}{gwire}{gdash}
    \pgoutgoing{0/gwire/{X_2},1/gdash/{X_2^*},2/gwire/{X_1},3/gdash/{X_1^*}}
    \node[glab] at (1.5,-1.4) {$c_{X\otimes X^*,\,X\otimes X^*}$};
  \end{scope}
\end{tikzpicture}
\caption{The morphism $S_2(c_{X,X})=s_2^{-1}(c_{X,X}\otimes\id_{X^*\otimes X^*})\,s_2$ and the braiding of $X\otimes X^*$ with itself, read from bottom to top, with the conventions of Figure~\ref{fig:dual}. The subscripts identify the two pairs of strands at the bottom. On the left the two solid strands exchange their places and the dashed strands keep their order. On the right both pairs of strands exchange their places.}
\label{fig:gate}
\end{figure}

The morphism $S_2(c_{X,X})$ crosses the two copies of $X$ in $X\otimes X^*\otimes X\otimes X^*$ and preserves the order of the two copies of $X^*$. In general, $S_2(c_{X,X})$ need not coincide with the braiding $c_{X\otimes X^*,X\otimes X^*}$, whose braid diagram also exchanges the two copies of $X^*$ (Figure~\ref{fig:gate}).

\subsection{The homomorphisms \texorpdfstring{$\Phi_n$}{Phi\_n}}\label{subsec:sufficiency}

We now compose the homomorphisms $S_n$ with the homomorphisms of Corollary~\ref{cor:regular}.

\begin{theorem}\label{thm:sufficiency}
Suppose that $d^2\in\Z$, and put $m:=d^2\FPdim(\cA)$. Let $R_\cA$ be the regular object of $\cA$, let \mbox{$W:=\Hom_\cA(R_\cA,X\otimes X^*\otimes R_\cA)$}, and let $\Phi^{X\otimes X^*}_n\colon\End_\cA((X\otimes X^*)^{\otimes n})\to\End(W^{\otimes n})$ be the homomorphisms of Corollary~\textup{\ref{cor:regular}} for $\cE=\cA$ and $T=X\otimes X^*$. Put
\[
\Phi_n:=\Phi^{X\otimes X^*}_n\circ S_n\colon E_n\longrightarrow\End(W^{\otimes n}),\qquad n\geq1 .
\]
\begin{enumerate}[label=\textup{(\alph*)}]
\item The maps $\Phi_n$ are injective unital algebra homomorphisms, they are compatible with tensor products in the sense of \eqref{eq:tc}, and $\dim W=m$.
\item For $n\geq1$ and each simple constituent $Y$ of $X^{\otimes n}$, put $K_{n,Y}:=\Hom_\cA(R_\cA,Y\otimes(X^*)^{\otimes n}\otimes R_\cA^{\otimes n})$. The object $Y\otimes(X^*)^{\otimes n}$ belongs to $\cA$, because it is a direct summand of $X^{\otimes n}\otimes(X^*)^{\otimes n}\cong (X\otimes X^*)^{\otimes n}$. Then $\dim K_{n,Y}=(m/d)^n\,d_Y$, and there is an isomorphism of $E_n$-modules
\begin{equation}\label{eq:module-decomposition}
W^{\otimes n}\cong\bigoplus_Y\Hom_\cC(Y,X^{\otimes n})\otimes K_{n,Y},
\end{equation}
where $Y$ runs over the simple constituents of $X^{\otimes n}$ and $f\in E_n$ acts on the right-hand side by $f_Y\otimes\id_{K_{n,Y}}$. In particular, the multiplicity of the simple $E_n$-module $\Hom_\cC(Y,X^{\otimes n})$ in $W^{\otimes n}$ is $(m/d)^n\,d_Y$.
\item For every $f\in E_n$ we have $m^{-n}\tr\bigl(\Phi_n(f)\bigr)=d^{-n}\Tr(f)$, where $\Tr$ is the trace of the canonical spherical structure of $\cD$.
\item The pair $(W,R)$ with $R:=\Phi_2(c_{X,X})$ is a localization of $X$ of dimension $m$.
\end{enumerate}
\end{theorem}

\begin{proof}
By Lemma~\ref{lem:integrality}(b), $\cA$ is an integral fusion category, and hence Corollary~\ref{cor:regular} applies to $\cE=\cA$ and $T=X\otimes X^*$. As in Section~\ref{sec:regular}, we put $H:=\Hom_\cA(R_\cA,-)$.

(a) The maps $\Phi^{X\otimes X^*}_n$ and $S_n$ are injective unital algebra homomorphisms by Corollary~\ref{cor:regular}(a) and Proposition~\ref{prop:dual-strands}(a). For $f\in E_n$ and $g\in E_k$, Proposition~\ref{prop:dual-strands}(b) and Corollary~\ref{cor:regular}(b) give
\[
\Phi_{n+k}(f\otimes g)=\Phi^{X\otimes X^*}_{n+k}\bigl(S_n(f)\otimes S_k(g)\bigr)=\Phi^{X\otimes X^*}_n\bigl(S_n(f)\bigr)\otimes\Phi^{X\otimes X^*}_k\bigl(S_k(g)\bigr)=\Phi_n(f)\otimes\Phi_k(g).
\]
Finally, $\dim W=d_{X\otimes X^*}\FPdim(\cA)$ by Corollary~\ref{cor:regular}(a), and $d_{X\otimes X^*}=d^2$.

(b) \emph{Step 1. The action of $E_n$ becomes composition.} In this proof, $\id$ without a subscript denotes the identity morphism of $(X^*)^{\otimes n}\otimes R_\cA^{\otimes n}$. Recall that $\Phi^{X\otimes X^*}_n(g)=J_n^{-1}H(g\otimes\id_{R_\cA^{\otimes n}})J_n$, where $J_n$ is an isomorphism from $W^{\otimes n}$ onto $\Hom_\cA(R_\cA,(X\otimes X^*)^{\otimes n}\otimes R_\cA^{\otimes n})$. Composition with the isomorphism $s_n\otimes\id_{R_\cA^{\otimes n}}$ is an isomorphism from this space onto $\Hom_\cA(R_\cA,X^{\otimes n}\otimes(X^*)^{\otimes n}\otimes R_\cA^{\otimes n})$. Let
\[
\Theta_n\colon W^{\otimes n}\longrightarrow\Hom_\cA\bigl(R_\cA,X^{\otimes n}\otimes(X^*)^{\otimes n}\otimes R_\cA^{\otimes n}\bigr),\qquad \Theta_n(w):=(s_n\otimes\id_{R_\cA^{\otimes n}})\,J_n(w),
\]
be the composite isomorphism. Since $S_n(f)=s_n^{-1}(f\otimes\id_{(X^*)^{\otimes n}})s_n$, the definition of $\Phi_n$ gives
\begin{equation}\label{eq:theta}
\Theta_n\bigl(\Phi_n(f)w\bigr)=(f\otimes\id)\,\Theta_n(w)\qquad(f\in E_n,\ w\in W^{\otimes n}).
\end{equation}
Thus $\Theta_n$ is an isomorphism of $E_n$-modules, where $f\in E_n$ acts on the target of $\Theta_n$ by composition with $f\otimes\id$.

\emph{Step 2. Decomposition of $X^{\otimes n}$.} For a simple constituent $Y$ of $X^{\otimes n}$, a morphism $u\colon Y\to X^{\otimes n}$ and an element $k\in K_{n,Y}$, the composite $(u\otimes\id)\,k$ is a morphism from $R_\cA$ to $X^{\otimes n}\otimes(X^*)^{\otimes n}\otimes R_\cA^{\otimes n}$. This defines a linear map
\[
\bigoplus_Y\Hom_\cC(Y,X^{\otimes n})\otimes K_{n,Y}\longrightarrow\Hom_\cA\bigl(R_\cA,X^{\otimes n}\otimes(X^*)^{\otimes n}\otimes R_\cA^{\otimes n}\bigr),\qquad u\otimes k\longmapsto(u\otimes\id)\,k .
\]
It is an isomorphism. Indeed, the same formula defines such a map for every direct summand of $X^{\otimes n}$ in place of $X^{\otimes n}$, these maps are compatible with direct sums, and for a simple summand $Y$ the map is the identity of $K_{n,Y}$, because $\Hom_\cC(Y,Y)=\C\,\id_Y$ and $\Hom_\cC(Y',Y)=0$ for $Y'\neq Y$. It is a homomorphism of $E_n$-modules, because $(f\otimes\id)(u\otimes\id)k=(fu\otimes\id)k$ and $fu=f_Y(u)$. In Figure~\ref{fig:module}, $f$ acts on the strand $X^{\otimes n}$ and meets only the box $u$. Together with Step~1 this gives \eqref{eq:module-decomposition}.

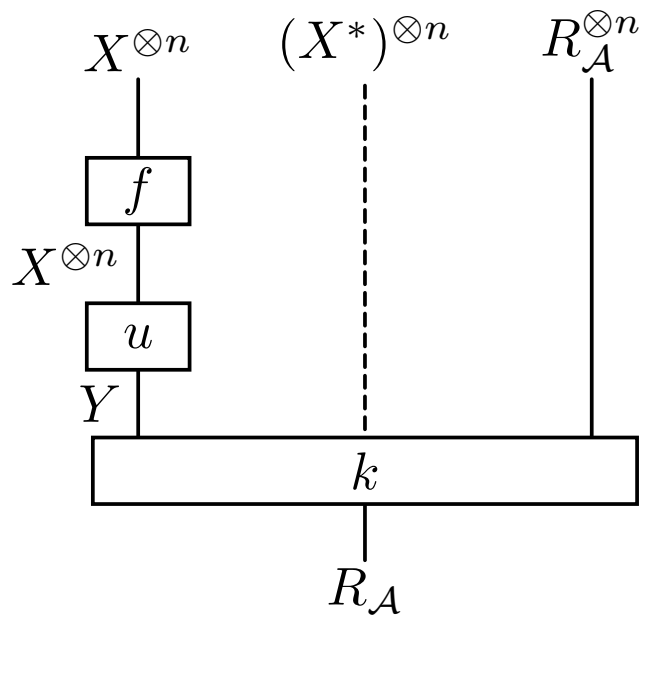
\begin{figure}[H]
\centering
\begin{tikzpicture}[x=1cm,y=.74cm]
  \draw[gwire] (1.5,0) -- (1.5,.8);
  \draw[gwire] (0,.8) -- (0,4.3);
  \draw[gdash] (1.5,.8) -- (1.5,4.3);
  \draw[gwire] (3,.8) -- (3,4.3);
  \node[gbox,minimum width=3.6cm] at (1.5,.8) {$k$};
  \node[gbox] at (0,2.0) {$u$};
  \node[gbox] at (0,3.3) {$f$};
  \node[glab,below] at (1.5,0) {$R_\cA$};
  \node[glab,anchor=east] at (-.06,1.4) {$Y$};
  \node[glab,anchor=east] at (-.06,2.65) {$X^{\otimes n}$};
  \node[glab,above] at (0,4.3) {$X^{\otimes n}$};
  \node[glab,above] at (1.5,4.3) {$(X^*)^{\otimes n}$};
  \node[glab,above] at (3,4.3) {$R_\cA^{\otimes n}$};
\end{tikzpicture}
\caption{The morphism $(fu\otimes\id)\,k$, for $k\in K_{n,Y}$, $u\colon Y\to X^{\otimes n}$ and $f\in E_n$.}
\label{fig:module}
\end{figure}

\emph{Step 3. The dimension of $K_{n,Y}$.} By Lemma~\ref{lem:beta}(a),
\[
\dim K_{n,Y}=d_{Y\otimes(X^*)^{\otimes n}\otimes R_\cA^{\otimes n}}=d_Y\,d^n\,\FPdim(\cA)^n=\Bigl(\frac md\Bigr)^nd_Y ,
\]
since $d_{R_\cA}=\FPdim(\cA)$ and $m/d=d\,\FPdim(\cA)$. The last assertion of (b) follows, because the spaces $\Hom_\cC(Y,X^{\otimes n})$ are the simple $E_n$-modules by \eqref{eq:wedderburn}.

(c) \emph{Step 1. The trace of $H(h)$.} Let $V$ be an object of $\cA$ and let $h\in\End_\cA(V)$. We claim that
\begin{equation}\label{eq:trace-H}
\tr\bigl(H(h)\bigr)=\Tr(h).
\end{equation}
Indeed, $R_\cA$ is the direct sum of $d_i$ copies of each simple object $i$ of $\cA$. Hence $H(V)$ is the direct sum of $d_i$ copies of each space $\Hom_\cA(i,V)$, and $H(h)$ acts on each of these copies by $h_i$. Therefore $\tr(H(h))=\sum_id_i\tr(h_i)$. This number is $\Tr(h)$ by Lemma~\ref{lem:trace-formula} in the fusion category $\cD$ with its canonical spherical structure, for which $\dim(i)=d_i$ by Lemma~\ref{lem:integrality}(b). The simple objects of $\cD$ that do not belong to $\cA$ do not contribute, because they are not constituents of $V$.

\emph{Step 2. Closing the additional strands.} By \eqref{eq:theta}, the isomorphism $\Theta_n$ identifies the operator $\Phi_n(f)$ with $H(f\otimes\id)$, and $X^{\otimes n}\otimes(X^*)^{\otimes n}\otimes R_\cA^{\otimes n}$ is an object of $\cA$. Hence \eqref{eq:trace-H} and the multiplicativity of the trace \cite[Proposition~4.7.3]{EGNO15} give
\[
\tr\bigl(\Phi_n(f)\bigr)=\Tr(f\otimes\id)=\Tr(f)\,\dim((X^*)^{\otimes n})\,\dim(R_\cA^{\otimes n})=\Tr(f)\,d^n\,\FPdim(\cA)^n=\Bigl(\frac md\Bigr)^n\Tr(f)
\]
(Figure~\ref{fig:closing}). Here we used that $\dim(V)=d_V$ for every object $V$ of $\cD$, that $d_{X^*}=d$ and that $d_{R_\cA}=\FPdim(\cA)$. The identity (c) also follows from (b) and Lemma~\ref{lem:trace-formula}, since $\dim K_{n,Y}=(m/d)^nd_Y$.

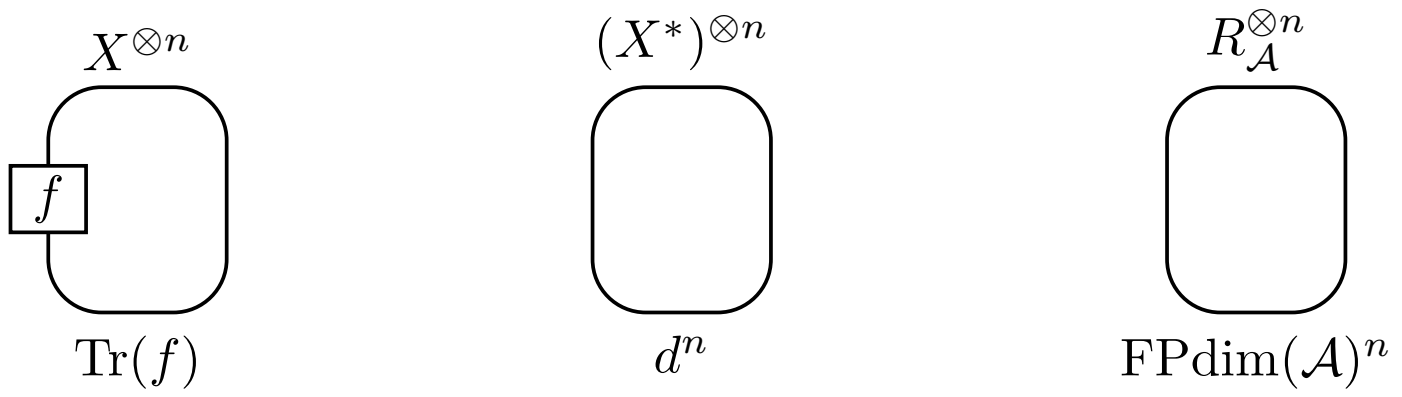
\begin{figure}[H]
\centering
\begin{tikzpicture}[x=1cm,y=1cm]
\foreach \x/\obj/\val in {0/{X^{\otimes n}}/{\Tr(f)},3.6/{(X^*)^{\otimes n}}/{d^n},7.4/{R_{\cA}^{\otimes n}}/{\FPdim(\cA)^n}}{
 \begin{scope}[xshift=\x cm]
 \draw[gwire,rounded corners=10pt] (0,.72)--(0,1.49)--(1.18,1.49)--(1.18,0)--(0,0)--cycle;
 \node[glab,above] at (.59,1.55) {$\obj$};
 \node[glab,below] at (.59,-.1) {$\val$};
 \end{scope}}
\node[gbox,minimum width=.5cm] at (0,.75) {$f$};
\end{tikzpicture}
\caption{The trace of $f\otimes\id$, where $\id$ is the identity morphism of $(X^*)^{\otimes n}\otimes R_\cA^{\otimes n}$. Closing the strands of $(X^*)^{\otimes n}$ gives the factor $d^n$, and closing the strands of $R_\cA^{\otimes n}$ gives the factor $\FPdim(\cA)^n$. Their product is $(m/d)^n$.}
\label{fig:closing}
\end{figure}

(d) This follows from (a) and Lemma~\ref{lem:localization}(b).
\end{proof}

\begin{figure}[H]
\centering
\begin{tikzpicture}[x=1.15cm,y=.76cm]
\draw[gwire] (0,0)--(0,4.95);
\draw[gwire] (2.2,0)--(2.2,1.03);
\draw[gwire] (2.2,3.92)--(2.2,4.95);
\foreach \x in {1.05,2.2,3.35}{\draw[gwire] (\x,1.40)--(\x,3.55);}
\node[gbox,minimum width=3.2cm,minimum height=.5cm] at (2.2,1.2) {$\widehat v$};
\node[gbox,minimum width=3.2cm,minimum height=.5cm] at (2.2,3.75) {$p_\lambda$};
\node[gbox,minimum width=1.9cm,minimum height=.5cm] at (.525,2.49) {$S_2(c_{X,X})$};
\node[glab,below] at (0,0) {$X\otimes X^*$}; \node[glab,above] at (0,4.95) {$X\otimes X^*$};
\node[glab,below] at (2.2,0) {$R_\cA$}; \node[glab,above] at (2.2,4.95) {$R_\cA$};
\node[glab,right,font=\scriptsize] at (1.05,1.78) {$X\otimes X^*$};
\node[glab,right,font=\scriptsize] at (2.2,2.5) {$R_\cA$};
\node[glab,right,font=\scriptsize] at (3.35,2.5) {$R_\cA$};
\end{tikzpicture}
\caption{The endomorphism of $X\otimes X^*\otimes R_\cA$ on the right-hand side of \eqref{eq:R-block}. The two strands $R_\cA$ that leave $\widehat v$ do not meet $S_2(c_{X,X})$. Composing the displayed endomorphism with $u$ gives the right-hand side of \eqref{eq:R-block}.}
\label{fig:block}
\end{figure}
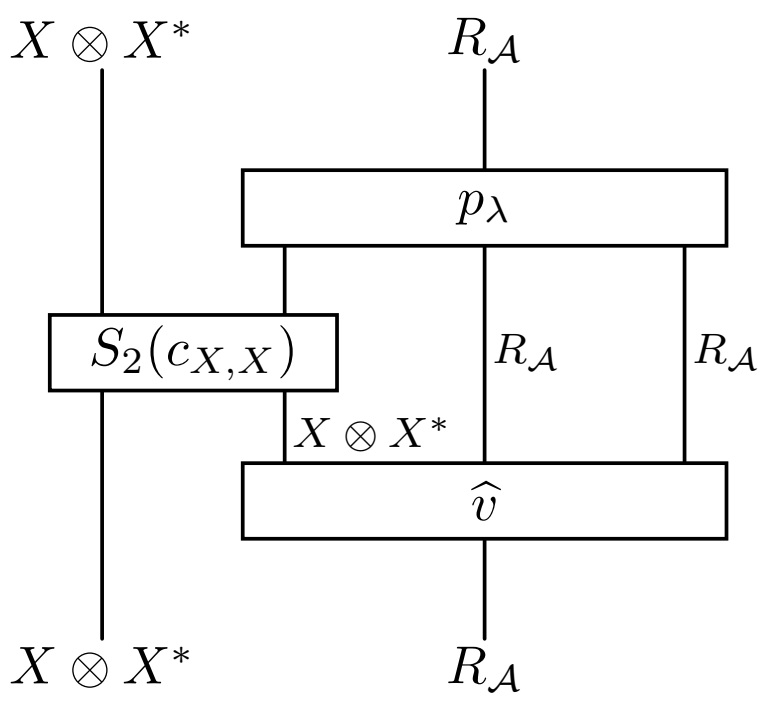

\begin{remark}\label{rem:explicit-R}
The operator $R$ of Theorem~\ref{thm:sufficiency} can be computed as follows. We have $S_2(c_{X,X})=s_2^{-1}(c_{X,X}\otimes\id_{X^*\otimes X^*})s_2$, where $s_2=\id_X\otimes c_{X^*,X}\otimes\id_{X^*}$ (Figure~\ref{fig:gate}). The isomorphism $J_2\colon W\otimes W\to\Hom_\cA(R_\cA,X\otimes X^*\otimes X\otimes X^*\otimes R_\cA\otimes R_\cA)$ is given by $J_2(u\otimes v)=(\id_{X\otimes X^*}\otimes\widehat v)\,u$, and $R$ is the linear map that corresponds, under $J_2$, to composition with $S_2(c_{X,X})\otimes\id_{R_\cA\otimes R_\cA}$.

The matrix of $R$ can be computed without inverting $J_2$. Let $\lambda$ be a linear form on $W$, and put $V:=X\otimes X^*\otimes R_\cA$. Then $W=H(V)$. By the Yoneda lemma, the maps $(\lambda\otimes\id)\,\Gamma_{S,V}^{-1}$ from $\Hom_\cA(S,V\otimes R_\cA)$ to $\Hom_\cA(S,R_\cA)$, which are natural in $S$, are given by composition with a unique morphism $p_\lambda\colon V\otimes R_\cA\to R_\cA$. Since $\widehat v=\Gamma_{R_\cA,V}(v\otimes\id_{R_\cA})$, we have $p_\lambda\,\widehat v=\lambda(v)\,\id_{R_\cA}$ for every $v\in W$. Then, for all $u,v\in W$,
\begin{equation}\label{eq:R-block}
(\id_W\otimes\lambda)\bigl(R(u\otimes v)\bigr)=(\id_{X\otimes X^*}\otimes p_\lambda)\,\bigl(S_2(c_{X,X})\otimes\id_{R_\cA\otimes R_\cA}\bigr)\,(\id_{X\otimes X^*}\otimes\widehat v)\,u
\end{equation}
(Figure~\ref{fig:block}). Indeed, $(\id_{X\otimes X^*}\otimes p_\lambda)\,J_2(x\otimes y)=\lambda(y)\,x$ for all $x,y\in W$, and $J_2\bigl(R(u\otimes v)\bigr)=\bigl(S_2(c_{X,X})\otimes\id_{R_\cA\otimes R_\cA}\bigr)J_2(u\otimes v)$. If $\cC$ is unitary, with the choices of Theorem~\ref{thm:unitary}, then for the linear form $\lambda=\langle v',-\rangle$ we have $p_\lambda=\widehat{v'}^{\,*}$ by \eqref{eq:delta-unitary}, and \eqref{eq:R-block} gives the matrix entries
\[
\bigl\langle u'\otimes v',R(u\otimes v)\bigr\rangle=\omega\Bigl(u'^{\,*}\,(\id_{X\otimes X^*}\otimes\widehat{v'}^{\,*})\,\bigl(S_2(c_{X,X})\otimes\id_{R_\cA\otimes R_\cA}\bigr)\,(\id_{X\otimes X^*}\otimes\widehat v)\,u\Bigr).
\]
\end{remark}

\begin{remark}\label{rem:dimension}
The dimension $m=d^2\FPdim(\cA)$ is not in general the smallest dimension of a localization. Suppose that $\cA$ admits a fiber functor $F\colon\cA\to\Vect$. Put $W_0:=F(X\otimes X^*)$, and let $\xi_n\colon W_0^{\otimes n}\to F((X\otimes X^*)^{\otimes n})$ be the iterated tensor structure of $F$. The maps $f\mapsto\xi_n^{-1}\,F(S_n(f))\,\xi_n$ are injective unital algebra homomorphisms $E_n\to\End(W_0^{\otimes n})$, because $F$ is faithful and $S_n$ is injective, and they are compatible with tensor products by Proposition~\ref{prop:dual-strands}(b) and the argument of Lemma~\ref{lem:transport}(a) below. By Lemma~\ref{lem:localization}, they give a localization of $X$ on $W_0$. Its dimension is $\dim F(X\otimes X^*)=d_{X\otimes X^*}=d^2$, since tensor functors preserve Frobenius--Perron dimensions \cite[Proposition~4.5.7]{EGNO15}. For the Ising anyon $\sigma$ we have $d^2=2$ and $\cA=\langle\psi\rangle$, with $\FPdim(\cA)=2$. Theorem~\ref{thm:sufficiency} gives a localization of dimension $4$. The fusion category $\langle\psi\rangle$ is braided equivalent to the fusion category of super vector spaces \cite[Lemma~B.11]{DGNO10}. Its underlying fusion category is the fusion category of $\Z/2$-graded vector spaces, and forgetting the grading is a fiber functor. The localization obtained from it has dimension $2$, the dimension of the localization of \cite{FRW06,RW12}.
\end{remark}

%======================================================================
\section{Unitary localizations}\label{sec:unitary}

We keep the notation of Section~\ref{sec:sufficiency}. The homomorphisms $\Phi_n$ of Theorem~\ref{thm:sufficiency} are composites of the homomorphisms $S_n$, defined with the dual strands, and the homomorphisms $\Phi^{X\otimes X^*}_n$, defined with the regular object. Both preserve adjoints when $\cC$ is unitary, and this gives the unitary case of Theorem~\ref{thm:main} (Theorem~\ref{thm:unitary}). We then show that it is enough to assume that $\cA$ is unitarizable (Corollary~\ref{cor:unitarizable}), and we prove the results stated in the introduction.

\begin{theorem}\label{thm:unitary}
Let $X$ be a nonzero object of a unitary braided fusion category $\cC$, suppose that $d^2\in\Z$, and let $\cA$ be as in Theorem~\textup{\ref{thm:sufficiency}}. Fix a faithful positive linear functional $\omega$ on $\End_\cA(R_\cA)$ and a morphism $\delta$ as in Lemma~\textup{\ref{lem:beta}(c)}, for $\cE=\cA$. Then $W=\Hom_\cA(R_\cA,X\otimes X^*\otimes R_\cA)$ is a Hilbert space with the inner product $\langle h,h'\rangle=\omega(h^*h')$, and the homomorphisms $\Phi_n$ of Theorem~\textup{\ref{thm:sufficiency}} satisfy $\Phi_n(f^*)=\Phi_n(f)^*$ for all $n\geq1$ and $f\in E_n$. In particular, $R=\Phi_2(c_{X,X})$ is unitary, and $(W,R)$ is a unitary localization of $X$ of dimension $d^2\FPdim(\cA)$.
\end{theorem}

\begin{proof}
The fusion subcategory $\cA$ of $\cC$ is unitary, with the adjoints of $\cC$, and it is integral by Lemma~\ref{lem:integrality}(b). By Lemma~\ref{lem:beta}(c), $W$ is a Hilbert space, and by Corollary~\ref{cor:regular}(c), $\Phi^{X\otimes X^*}_n(g^*)=\Phi^{X\otimes X^*}_n(g)^*$ for every $g\in\End_\cA((X\otimes X^*)^{\otimes n})$. By Proposition~\ref{prop:dual-strands}(c), $S_n(f^*)=S_n(f)^*$. Hence
\[
\Phi_n(f^*)=\Phi^{X\otimes X^*}_n\bigl(S_n(f)^*\bigr)=\Phi^{X\otimes X^*}_n\bigl(S_n(f)\bigr)^*=\Phi_n(f)^* .
\]
The operator $R$ is unitary by Lemma~\ref{lem:localization}(c), and $(W,R)$ is a localization of dimension $d^2\FPdim(\cA)$ by Theorem~\ref{thm:sufficiency}.
\end{proof}

\subsection{Transport along a tensor equivalence}\label{subsec:transport}
In Theorem~\ref{thm:unitary} only the unitary structure of the fusion category $\cD=\langle X\rangle$ was used, since every braiding of a unitary fusion category is unitary. In this subsection $\cC$ is a braided fusion category that is not necessarily unitary. We show that a weakly integral object $X$ admits a unitary localization whenever $\cD$ is unitarizable (Proposition~\ref{prop:transport}), and that $\cD$ is unitarizable if and only if $\cA$ is unitarizable (Proposition~\ref{prop:unitarizable}). The method is to transport the braiding of $\cD$ to the unitary fusion category and to apply Theorem~\ref{thm:unitary} there.

Let $G\colon\cD\to\cD'$ be a tensor equivalence from $\cD$ to a fusion category $\cD'$, with tensor structure $\gamma_{A,B}\colon G(A)\otimes G(B)\to G(A\otimes B)$. Put $X':=G(X)$. We give $\cD'$ the braiding $c'$ transported along $G$, which is determined by
\begin{equation}\label{eq:transported-braiding}
c'_{G(A),G(B)}=\gamma_{B,A}^{-1}\,G(c_{A,B})\,\gamma_{A,B}.
\end{equation}
As in Subsection~\ref{subsec:braided}, we write $\gamma_{n,k}:=\gamma_{X^{\otimes n},X^{\otimes k}}$ and $\gamma_n\colon X'^{\otimes n}\to G(X^{\otimes n})$ for the iterated structure isomorphisms of $G$ at the object $X$. Since $G$ is an equivalence, the maps
\begin{equation}\label{eq:Gn}
G_n\colon E_n\longrightarrow\End_{\cD'}(X'^{\otimes n}),\qquad
G_n(f):=\gamma_n^{-1}\,G(f)\,\gamma_n
\end{equation}
are algebra isomorphisms.

\begin{lemma}\label{lem:transport}
Let $G\colon\cD\to\cD'$ be a tensor equivalence and let $G_n$ be the isomorphisms \eqref{eq:Gn}.
\begin{enumerate}[label=\textup{(\alph*)}]
\item For $f\in E_n$ and $g\in E_k$, we have $G_{n+k}(f\otimes g)=G_n(f)\otimes G_k(g)$.
\item For every $n\geq2$, we have $G_n\circ\rho_n^X=\rho_n^{X'}$, where $\rho^{X'}_n$ is defined with the braiding $c'$.
\end{enumerate}
\end{lemma}

\begin{proof}
(a) By the naturality of $\gamma$, we have $G(f\otimes g)\,\gamma_{n,k}=\gamma_{n,k}\,(G(f)\otimes G(g))$. By \eqref{eq:gamma-iterated} for the object $X$ and this identity,
\[
G_{n+k}(f\otimes g)=\bigl(\gamma_n\otimes\gamma_k\bigr)^{-1}\bigl(G(f)\otimes G(g)\bigr)\bigl(\gamma_n\otimes\gamma_k\bigr)=G_n(f)\otimes G_k(g).
\]

(b) By \eqref{eq:transported-braiding} for $A=B=X$, we have $G_2(c_{X,X})=c'_{X',X'}$. Since $\rho^X_n(\sigma_i)$ is the tensor product of $c_{X,X}$ with identity morphisms and the homomorphisms $G_k$ are unital, (a) gives $G_n(\rho^X_n(\sigma_i))=\rho^{X'}_n(\sigma_i)$ for every generator. Both sides of (b) are algebra homomorphisms, and hence they are equal.
\end{proof}

\begin{proposition}\label{prop:transport}
Let $X$ be a nonzero object of a braided fusion category $\cC$, suppose that $d^2\in\Z$, and let $G\colon\cD\to\cD'$ be a tensor equivalence onto a unitary fusion category, with the isomorphisms $G_n$ of \eqref{eq:Gn}. Then there are a Hilbert space $W$ of dimension $d^2\FPdim(\cA)$ and injective unital algebra homomorphisms $\Phi_n\colon E_n\to\End(W^{\otimes n})$, $n\geq1$, that are compatible with tensor products and such that the homomorphisms $\Phi_n\circ G_n^{-1}\colon\End_{\cD'}(X'^{\otimes n})\to\End(W^{\otimes n})$ preserve adjoints. The operator $R:=\Phi_2(c_{X,X})$ is unitary, and $(W,R)$ is a unitary localization of~$X$.
\end{proposition}

\begin{proof}
\emph{Step 1. The fusion category $\cD'$.} The isomorphisms \eqref{eq:transported-braiding} define a braiding $c'$ on $\cD'$, because $G$ is a tensor equivalence. Thus $\cD'$ is a unitary braided fusion category. A tensor equivalence preserves duals and Frobenius--Perron dimensions. Hence $d_{X'}=d$, and $G(X^{\otimes n})\cong X'^{\otimes n}$ for every $n\geq1$. Since $G$ is essentially surjective, every simple object of $\cD'$ is isomorphic to $G(Y)$ for a simple object $Y$ of $\cD$, and $Y$ is a direct summand of some $X^{\otimes n}$ by Lemma~\ref{lem:generated}. Hence every simple object of $\cD'$ is a direct summand of some $X'^{\otimes n}$, and $\cD'=\langle X'\rangle$ by Lemma~\ref{lem:generated}. Put $\cA':=\langle X'\otimes X'^*\rangle$. Then $X'\otimes X'^*\cong G(X\otimes X^*)$, and by Lemma~\ref{lem:generated}, applied to $X\otimes X^*$ and to $X'\otimes X'^*$, the functor $G$ sends the simple objects of $\cA$ to simple objects of $\cA'$, and every simple object of $\cA'$ is isomorphic to the image of a simple object of $\cA$. Since $G$ preserves dimensions, $\FPdim(\cA')=\FPdim(\cA)$.

\emph{Step 2. The homomorphisms.} Theorem~\ref{thm:unitary}, applied to the object $X'$ of $\cD'$, gives a Hilbert space $W$ of dimension $d^2\FPdim(\cA)$ and injective unital algebra homomorphisms $\Phi'_n\colon\End_{\cD'}(X'^{\otimes n})\to\End(W^{\otimes n})$ that are compatible with tensor products and preserve adjoints. Put $\Phi_n:=\Phi'_n\circ G_n$. These maps are injective unital algebra homomorphisms, and they are compatible with tensor products by Lemma~\ref{lem:transport}(a). Moreover, $\Phi_n\circ G_n^{-1}=\Phi'_n$ preserves adjoints.

\emph{Step 3. The localization.} By Lemma~\ref{lem:localization}(b), $(W,R)$ is a localization of $X$. Finally, $R=\Phi'_2(G_2(c_{X,X}))=\Phi'_2(c'_{X',X'})$ by Lemma~\ref{lem:transport}(b), and this operator is unitary by Theorem~\ref{thm:unitary}.
\end{proof}

We now compare the unitarizability of $\cD$ and of $\cA$. Recall that the \emph{adjoint subcategory} $\cD_\ad$ of $\cD$ is the fusion subcategory generated by the objects $Y\otimes Y^*$, where $Y$ runs over the simple objects of $\cD$. It is the neutral component of the universal grading of $\cD$, which is faithful \cite[Theorem~3.5 and Corollary~3.7]{GN08}. By Lemma~\ref{lem:integrality}(a), $\cD_\ad\subseteq\cA$. If $X$ is simple, then $\cA=\cD_\ad$, because $X\otimes X^*$ is one of the generators of $\cD_\ad$. For an object that is not simple the inclusion can be proper. For example, let $X=\one\oplus g$ in the fusion category of $\Z/2$-graded vector spaces, where $g$ is the simple object of degree one. Then $\cD_\ad$ consists of the multiples of $\one$, and $\cA=\cD$.

We use the following notion from \cite[Definition~5.13]{GHR13}. A fusion category $\mathcal K$ is \emph{completely unitary} if it has the following three properties.
\begin{enumerate}[label=\textup{(\roman*)}]
\item The fusion category $\mathcal K$ is tensor equivalent to a unitary fusion category, which is unique up to unitary tensor equivalence.
\item Every module category over $\mathcal K$ is equivalent to a unitary module category, which is unique up to unitary module equivalence.
\item Every module equivalence between unitary module categories over $\mathcal K$ is isomorphic to a unitary module equivalence, which is unique up to unitary module natural isomorphism.
\end{enumerate}
Here the module categories are finite and semisimple, and in \textup{(ii)} and \textup{(iii)} they are module categories over the unitary fusion category of \textup{(i)}.

\begin{lemma}\label{lem:completely-unitary}
Every unitarizable fusion category is completely unitary.
\end{lemma}

\begin{proof}
Let $\mathcal K$ be a unitarizable fusion category. Property \textup{(i)} holds because the unitary structure of $\mathcal K$ is unique up to unitary tensor equivalence \cite[Theorem~1]{Reu23}. For property \textup{(ii)}, a finite semisimple module category over a unitary fusion category is a direct sum of indecomposable module categories. Each summand has a unitary structure \cite[Corollary~3.8]{CGGH23}, and the orthogonal direct sum of these structures is a unitary structure on the module category. This structure is unique up to unitary module equivalence \cite[Corollary~9]{Reu23}. The result of \cite{CGGH23} is stated for right module categories, and it applies to left module categories after reversing the tensor product. For property \textup{(iii)}, each module equivalence between unitary module categories is isomorphic to a unitary module equivalence, and two such choices differ by a unitary isomorphism, by polar decomposition \cite[Theorem~8 and Remark~10]{Reu23}.
\end{proof}

\begin{proposition}\label{prop:unitarizable}
For every nonzero object $X$ of a braided fusion category, the fusion category $\langle X\rangle$ is unitarizable if and only if $\langle X\otimes X^*\rangle$ is unitarizable.
\end{proposition}

\begin{proof}
A fusion subcategory of a unitarizable fusion category is unitarizable, and $\cA$ is a fusion subcategory of $\cD$. Hence it suffices to prove that $\cD$ is unitarizable if $\cA$ is. Suppose that $\cA$ is unitarizable.

\emph{Step 1. The fusion category $\cD_\ad\boxtimes\cD_\ad^{\rm rev}$ is unitarizable.} The fusion subcategory $\cD_\ad$ of $\cA$ is unitarizable. Let $\cD_\ad^{\rm rev}$ denote $\cD_\ad$ with the reversed tensor product. If $\mathcal K'$ is a unitary fusion category tensor equivalent to $\cD_\ad$, then $\mathcal K'\boxtimes(\mathcal K')^{\rm rev}$ is a unitary fusion category tensor equivalent to $\cD_\ad\boxtimes\cD_\ad^{\rm rev}$. Hence $\cD_\ad\boxtimes\cD_\ad^{\rm rev}$ is unitarizable.

\emph{Step 2. Complete unitarity.} By Lemma~\ref{lem:completely-unitary}, the fusion categories $\cD_\ad$ and $\cD_\ad\boxtimes\cD_\ad^{\rm rev}$ are completely unitary.

\emph{Step 3. The extension theorem.} The universal grading is faithful, and it presents $\cD$ as a faithfully graded fusion category with neutral component $\cD_\ad$. By the extension theorem \cite[Theorem~5.17]{GHR13}, a faithfully graded fusion category is completely unitary when its neutral component $\mathcal K$ and the fusion category $\mathcal K\boxtimes\mathcal K^{\rm rev}$ are completely unitary. By Step~2, $\cD$ is completely unitary, and in particular it is unitarizable.
\end{proof}

\begin{corollary}\label{cor:unitarizable}
Let $X$ be a nonzero weakly integral object of a braided fusion category such that $\cA=\langle X\otimes X^*\rangle$ is unitarizable. Then $X$ admits a unitary localization of dimension $d^2\FPdim(\cA)$.
\end{corollary}

\begin{proof}
By Proposition~\ref{prop:unitarizable}, the fusion category $\cD=\langle X\rangle$ is unitarizable. Let $G\colon\cD\to\cD'$ be a tensor equivalence onto a unitary fusion category. Proposition~\ref{prop:transport}, applied to $G$, gives a unitary localization of $X$ of dimension $d^2\FPdim(\cA)$.
\end{proof}

\subsection{Proofs of the results stated in the introduction}

\begin{proof}[Proof of Theorem~\textup{\ref{thm:main}}]
The number $m=d^2\FPdim(\cA)$ is an integer because $\cA$ is integral (Lemma~\ref{lem:integrality}(b)). Parts~(a) and~(b) of Theorem~\ref{thm:main} are Theorem~\ref{thm:sufficiency}(a)--(c), and part~(c) of Theorem~\ref{thm:main} is Theorem~\ref{thm:unitary}.
\end{proof}

\begin{proof}[Proof of Theorem~\textup{\ref{thm:rw}}]
The pair $(W,R)$ of Theorem~\ref{thm:sufficiency}(d) is a localization of $X$. If $\cC$ is unitary, it is a unitary localization for the choices of Theorem~\ref{thm:unitary}.
\end{proof}

\begin{proof}[Proof of Corollary~\textup{\ref{cor:rw-intro}}]
(a) Every simple object of a weakly integral fusion category is weakly integral \cite[Proposition~9.6.9(i)]{EGNO15}, and Theorem~\ref{thm:rw} applies.

(b) Let $X$ be a simple object such that $\cD=\langle X\rangle$ is weakly group-theoretical. Then $\cD$ is weakly integral \cite[Remark~9.8.2]{EGNO15}, and hence $d^2$ is an integer by Lemma~\ref{lem:integrality}(c). Moreover, $\cD$ is completely unitary \cite[Theorem~5.20]{GHR13}, and in particular it is tensor equivalent to a unitary fusion category. Proposition~\ref{prop:transport}, applied to such an equivalence, gives a unitary localization of $X$. For the last assertion, recall that every fusion subcategory of a weakly group-theoretical fusion category is weakly group-theoretical \cite[Proposition~9.8.4]{EGNO15}.
\end{proof}

\begin{remark}\label{rem:finite-order}
Let $X$ be a nonzero weakly integral object of a braided fusion category, and let $R=\Phi_2(c_{X,X})$ be the operator of Theorem~\ref{thm:sufficiency}. Then $R$ has finite order. Indeed, $\cD$ is pseudo-unitary by Lemma~\ref{lem:integrality}(b), and we consider its twist $\theta$ as in Subsection~\ref{subsec:ribbon}. The twists of the simple objects of $\cD$ are roots of unity \cite[Corollary~8.18.2(i)]{EGNO15}, and hence there is an integer $h\geq1$ such that $\theta_A^h=\id_A$ for every object $A$ of $\cD$. By the defining identity of a twist \cite[(8.32)]{EGNO15}, we have $c_{X,X}^2=\theta_{X\otimes X}\,(\theta_X^{-1}\otimes\theta_X^{-1})$, and the two factors commute because $\theta$ is natural. Hence $c_{X,X}^{2h}=\id_{X\otimes X}$, and therefore $R^{2h}=\id_{W\otimes W}$.

Since $\rho^R_n=\Phi_n\circ\rho^X_n$ and $\Phi_n$ is injective, the groups $\rho^X_n(B_n)$ and $\rho^R_n(B_n)$ are isomorphic. Therefore, for the unitary localization of Theorem~\ref{thm:unitary}, the conjecture of Rowell and Wang that a unitary Yang--Baxter operator of finite order has finite braid group images modulo scalars \cite[Conjecture~3.1]{RW12} would imply that the group $\Gamma:=\rho^R_n(B_n)$ is finite modulo scalars. For this group, finiteness modulo scalars is equivalent to finiteness. Indeed, let $q$ be the order of $R$ and put $N:=m^n=\dim W^{\otimes n}$. The determinant of each generator $\rho^R_n(\sigma_i)$ is a $q$-th root of unity, and hence the determinant of every element of $\Gamma$ is also a $q$-th root of unity. A scalar operator $\zeta\,\id_{W^{\otimes n}}$ in $\Gamma$ has determinant $\zeta^N$, and therefore $\zeta^{Nq}=1$. Thus $\Gamma$ contains finitely many scalar operators, and it is finite if it is finite modulo scalars. Consequently, the conjecture of Rowell and Wang would imply that the groups $\rho^X_n(B_n)$ are finite, which is the statement predicted by the Property~F conjecture for weakly integral objects \cite{NR11}, \cite[Conjecture~6.2]{RW18}. When $\cD$ is weakly group-theoretical, these groups are known to be finite \cite[Theorem~7.2]{GN21}. In the opposite direction, a simple object of a pseudo-unitary braided fusion category whose braid group images are finite is weakly integral \cite[Corollary~1.2]{LR26}.
\end{remark}

\begin{remark}\label{rem:wgt}
For a nonzero object of integral dimension $d$ in a weakly group-theoretical braided fusion category that admits a fiber functor, a unitary localization of dimension $d$ was constructed in \cite[Corollary~5.25]{GHR13}. Proposition~\ref{prop:transport}, together with \cite[Theorem~5.20]{GHR13}, gives a unitary localization of dimension $d^2\FPdim(\cA)$ for every nonzero weakly integral object $X$ such that $\langle X\rangle$ is weakly group-theoretical, whether $X$ is simple or not. It requires neither a fiber functor nor that $d$ be an integer. Corollary~\ref{cor:rw-intro}(b) is the case of a simple object.
\end{remark}

%======================================================================
\section{Partial traces and the Turaev enhancement}\label{sec:turaev}

In this section we prove Theorem~\ref{thm:enhanced-intro}. The proof uses only three properties of the homomorphisms of Theorem~\ref{thm:main}, and we state the assumptions accordingly. Let $X$ be a simple object of a unitary braided fusion category $\cC$, put $d:=d_X$, and suppose that $d^2\in\Z$. Let $W$ be a nonzero finite-dimensional Hilbert space, put $m:=\dim W$ and $\tau_n:=m^{-n}\tr$, and let $\Phi_n\colon E_n\to\End(W^{\otimes n})$, $n\geq1$, be unital algebra homomorphisms of the algebras $E_n=\End_\cC(X^{\otimes n})$ that are compatible with tensor products in the sense of \eqref{eq:tc}, preserve adjoints, and satisfy the trace identity
\begin{equation}\label{eq:trace-identity-6}
\tau_n\bigl(\Phi_n(f)\bigr)=d^{-n}\Tr(f),\qquad f\in E_n,
\end{equation}
where $\Tr$ is the trace of the canonical spherical structure of $\langle X\rangle$. For this structure $\dim(X)=d$. The homomorphisms of Theorem~\ref{thm:main}(c) have these properties, with $m=d^2\FPdim(\cA)$ and $\cA=\langle X\otimes X^*\rangle$. Put $R:=\Phi_2(c_{X,X})$. By Lemma~\ref{lem:localization}, $R$ is a unitary Yang--Baxter operator and $\rho^R_n=\Phi_n\circ\rho^X_n$ for $n\geq2$. The group $B_1$ is trivial, and we let $\rho^X_1$ and $\rho^R_1$ be the trivial homomorphisms, with values $\id_X$ and $\id_W$. The trace identity therefore gives
\begin{equation}\label{eq:common-braid-trace}
 \tau_n\bigl(\rho_n^R(b)\bigr)
       =d^{-n}\Tr\bigl(\rho_n^X(b)\bigr),\qquad b\in B_n,\ n\geq1 .
\end{equation}
Thus $d^n\,\tau_n\bigl(\rho^R_n(b)\bigr)$ is the Reshetikhin--Turaev invariant of the closure of $b$ with its blackboard framing and with every component colored by $X$ (Subsection~\ref{subsec:ribbon}).

We use the following notation for partial traces. We identify $\End(W\otimes W)$ with $\End(W)\otimes\End(W)$, and we denote by $\tau_1\otimes\id$ and $\id\otimes\tau_1$ the linear maps $\End(W\otimes W)\to\End(W)$ such that
\[
(\tau_1\otimes\id)(x\otimes y)=\tau_1(x)\,y,\qquad(\id\otimes\tau_1)(x\otimes y)=\tau_1(y)\,x\qquad(x,y\in\End(W)).
\]
These are the normalized partial traces over the first and over the second factor. In these expressions $\id$ is the identity map of $\End(W)$, whereas $\id_W$ is the identity operator of $W$.

\begin{proof}[Proof of Theorem~\textup{\ref{thm:enhanced-intro}}]
\emph{Step 1. Adding a kink.} Let $s\colon B_n\to B_{n+1}$ be the homomorphism with $s(\sigma_i)=\sigma_{i+1}$, which adds a strand on the left. We prove that, for every $n\geq1$, every $b\in B_n$ and each sign,
\begin{equation}\label{eq:stabilization}
\tau_{n+1}\bigl(\rho^R_{n+1}(\sigma_1^{\pm1}\,s(b))\bigr)=\frac{\theta_X^{\pm1}}{d}\,\tau_n\bigl(\rho^R_n(b)\bigr).
\end{equation}
The closure of $\sigma_1^{\pm1}s(b)$ is obtained from the closure of $b$ by adding a kink to the component that contains the first strand of $b$ (Figure~\ref{fig:stab}). Hence the two closures are isotopic as oriented links, and the blackboard framing of that component changes by $\pm1$. By Subsection~\ref{subsec:ribbon},
\[
\Tr\bigl(\rho^X_{n+1}(\sigma_1^{\pm1}\,s(b))\bigr)=\theta_X^{\pm1}\,\Tr\bigl(\rho^X_n(b)\bigr),
\]
and \eqref{eq:common-braid-trace}, for $n+1$ and for $n$, gives \eqref{eq:stabilization}. Figure~\ref{fig:stab} shows the same computation in terms of morphisms. The kink is the endomorphism of $X$ obtained from $c_{X,X}^{\pm1}$ by closing its first strand. It is a scalar because $X$ is simple, and the scalar is $\theta_X^{\pm1}$, since closing the second strand gives $\Tr(c_{X,X}^{\pm1})=\theta_X^{\pm1}\dim(X)$ by \eqref{eq:twist-trace}.

\begin{figure}[H]
\centering
\begin{tikzpicture}[x=1cm,y=.86cm]
 \begin{scope}
 \draw[gwire] (1.55,0) -- (1.55,1.31);
 \draw[gwire] (2.9,0) -- (2.9,2.95);
 \draw[gwire] (.55,2.31) -- (.55,2.47)
      .. controls (.55,2.90) and (-.11,2.90) .. (-.11,2.44)
      -- (-.11,1.18) .. controls (-.11,.76) and (.55,.76) .. (.55,1.16)
      -- (.55,1.31);
 \draw[gwire] (1.55,1.31) -- (.55,2.31);
 \draw[gwire] (.55,1.31) -- (1.55,2.31) -- (1.55,2.95);
 \draw[gwire,preaction={draw=white,line width=4.4pt}] (.87,1.63) -- (1.23,1.99);
 \node[gbox,minimum width=1.85cm] at (2.225,.54) {$b$};
 \node[glab,above] at (1.55,2.95) {$X$}; \node[glab,below] at (1.55,0) {$X$};
 \node[glab,above] at (2.9,2.95) {$X^{\otimes(n-1)}$}; \node[glab,below] at (2.9,0) {$X^{\otimes(n-1)}$};
 \end{scope}
 \node at (4.1,1.58) {$=$};
 \begin{scope}[xshift=5.0cm]
 \draw[gwire] (0,0) -- (0,2.95);
 \draw[gwire] (1.35,0) -- (1.35,2.95);
 \node[gbox,minimum width=1.85cm] at (.675,.54) {$b$};
 \node[gbox,minimum width=.72cm] at (0,1.83) {$\theta_X$};
 \node[glab,above] at (0,2.95) {$X$}; \node[glab,below] at (0,0) {$X$};
 \node[glab,above] at (1.35,2.95) {$X^{\otimes(n-1)}$}; \node[glab,below] at (1.35,0) {$X^{\otimes(n-1)}$};
 \end{scope}
\end{tikzpicture}
\caption{The braid $\sigma_1\,s(b)$ with its first strand closed, read from bottom to top. The box $b$ stands for $\rho^X_n(b)$. Closing the added strand leaves a kink on the first strand of $b$, which is the multiplication by $\theta_X$. Closing the remaining strands gives $\Tr\bigl(\rho^X_{n+1}(\sigma_1s(b))\bigr)=\theta_X\Tr\bigl(\rho^X_n(b)\bigr)$.}
\label{fig:stab}
\end{figure}
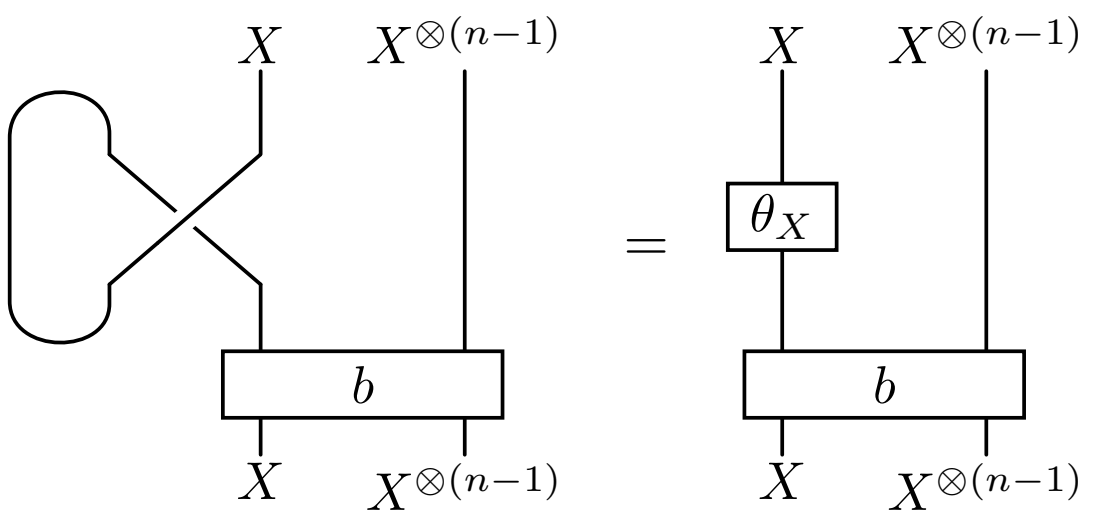

\emph{Step 2. The partial traces of $R$.} For $n=1$ and the trivial braid, \eqref{eq:stabilization} with the positive sign gives
\[
z:=\tau_2(R)=\frac{\theta_X}{d}.
\]
Let $B_n\subset B_{n+1}$ be the inclusion that adds a strand on the right. For $b\in B_n$ we have $\rho^R_{n+1}(b)=\rho^R_n(b)\otimes\id_W$, and hence $\tau_{n+1}\bigl(\rho^R_{n+1}(b)\bigr)=\tau_n\bigl(\rho^R_n(b)\bigr)$. Therefore the normalized traces define a function on the union $B_\infty$ of the groups $B_n$. We denote by $\tau_R\colon\C B_\infty\to\C$ its linear extension to the group algebra, which satisfies $\tau_R(xy)=\tau_R(yx)$. The homomorphisms $s$ are compatible with the inclusions and define a homomorphism $s\colon B_\infty\to B_\infty$, which we extend linearly to $\C B_\infty$. Since $\tau_R(\sigma_1)=z$, the identity \eqref{eq:stabilization} with the positive sign says that
\[
\tau_R\bigl(\sigma_1\,s(b)\bigr)=\tau_R(\sigma_1)\,\tau_R(b)
\]
for every $b\in B_\infty$, and by linearity for every $b\in\C B_\infty$. Thus $\tau_R$ is a Markov trace in the sense of \cite[Definition~2.2]{Lec26}. By a proposition of Lechner \cite[Proposition~2.3]{Lec26}, for a unitary Yang--Baxter operator $R$ the trace $\tau_R$ is a Markov trace if and only if $(\tau_1\otimes\id)(R)=\tau_2(R)\,\id_W$. Hence
\[
(\tau_1\otimes\id)(R)=z\,\id_W .
\]
The two normalized partial traces of a unitary Yang--Baxter operator coincide, by a theorem of Conti and Lechner \cite[Theorem~5.10(a)]{CL21}. Therefore $(\id\otimes\tau_1)(R)=z\,\id_W$ as well.

It remains to consider $R^{-1}$. Both partial traces commute with adjoints, that is, $(\tau_1\otimes\id)(a^*)=\bigl((\tau_1\otimes\id)(a)\bigr)^*$ and $(\id\otimes\tau_1)(a^*)=\bigl((\id\otimes\tau_1)(a)\bigr)^*$ for every $a\in\End(W\otimes W)$. Indeed, both sides are conjugate linear in $a$, and for $a=x\otimes y$ they are equal to $\overline{\tau_1(x)}\,y^*$ and to $\overline{\tau_1(y)}\,x^*$, respectively. Since $R^{-1}=R^*$, both partial traces of $R^{-1}$ are equal to $\overline z\,\id_W=\theta_X^{-1}d^{-1}\id_W$, because $\theta_X$ is a root of unity. Finally, $|\theta_X|=1$ gives $d=|\tau_2(R)|^{-1}$ and $\theta_X=d\,\tau_2(R)$. This proves (a).

\emph{Step 3. The enhancement.} In Turaev's convention \cite[Section~2.3]{Tur88}, an enhancement of $R$ consists of an operator $\mu\in\End(W)$ and scalars $\alpha,\beta\in\C^\times$ such that $\mu\otimes\mu$ commutes with $R$ and
\[
 (\id\otimes\tr)\bigl(R^{\pm1}(\mu\otimes\mu)\bigr)
       =\alpha^{\pm1}\beta\,\mu ,
\]
where $\id\otimes\tr$ is the partial trace over the second factor, without normalization. We take $\mu:=(d/m)\id_W$. Then $\mu\otimes\mu$ is a scalar operator and commutes with $R$. Since $\id\otimes\tr=m\,(\id\otimes\tau_1)$, part~(a) of Theorem~\ref{thm:enhanced-intro} gives
\[
 (\id\otimes\tr)\bigl(R^{\pm1}(\mu\otimes\mu)\bigr)
    =\Bigl(\frac dm\Bigr)^2m\,\frac{\theta_X^{\pm1}}d\,\id_W
    =\theta_X^{\pm1}\mu .
\]
Hence $(R,\mu,\alpha,\beta)=(R,(d/m)\id_W,\theta_X,1)$ is an enhanced Yang--Baxter operator.

\emph{Step 4. The knot invariant.} For a braid $b\in B_n$ with $n\geq2$, Turaev's invariant associated with this enhancement is $\alpha^{-w(b)}\beta^{-n}\tr(\rho_n^R(b)\mu^{\otimes n})$ \cite[Theorem~3.1.2 and Section~3.2]{Tur88}. For our enhancement it equals
\[
 \theta_X^{-w(b)}\Bigl(\frac dm\Bigr)^n\tr\bigl(\rho^R_n(b)\bigr)
 =\theta_X^{-w(b)}d^n\,\tau_n\bigl(\rho^R_n(b)\bigr)
   =\theta_X^{-w(b)}\Tr\bigl(\rho_n^X(b)\bigr),
\]
by \eqref{eq:common-braid-trace}. When the closure of $b$ is a knot, the right-hand side is $\RT_X(\widehat b)$ by \eqref{eq:rt}. This proves \eqref{eq:local-rt} for $n\geq2$. For $n=1$ the braid is trivial, it acts on $W$ by the identity, and its closure is the unknot. Both sides of \eqref{eq:local-rt} are then equal to $d$, since $\tr(\id_W)=m$.
\end{proof}

\section*{Acknowledgments}
The author was partially supported by Grant INV-2025-213-3452 from the School of Science of Universidad de los Andes. The author used Claude Fable~5.2 and Gemini~4 Carbon during the development and drafting of this paper, to explore and implement several proof strategies. All statements and proofs were checked by the author, who is responsible for them. The author thanks Anthropic for access to Claude through the Claude Team plan for scientists, and Google for access to Google Antigravity through its program for scientists.

\bibliographystyle{amsalpha}
\raggedbottom
\bibliography{weakly_integral_localization}

\end{document}